\documentclass[11pt]{amsart}
 \usepackage{graphicx}
\usepackage{xstring}
\usepackage{forloop}
\usepackage{bbm, amsmath, amsfonts, amscd, latexsym, amsthm, amssymb, graphicx, mathrsfs, comment, changepage}
\usepackage[colorlinks=true]{hyperref}

\usepackage{abstract}
\usepackage{mathtools}
\usepackage{caption}
\usepackage{MnSymbol}
\usepackage[dvipsnames]{xcolor}

\usepackage{tikz-cd}
\tikzset{
  symbol/.style={
    draw=none,
    every to/.append style={
      edge node={node [sloped, allow upside down, auto=false]{$#1$}}}
  }
}

\makeatletter
\newif\if@check@engine  \@check@enginetrue 
\makeatother
\usepackage[usenames,dvipsnames]{pstricks}

\usepackage{booktabs}

\newtheorem{theor}{\hspace{1cm}{\sc Theorem}}[section]
\newtheorem{utver}[theor]{\hspace{1cm}{\sc Proposition}}
\newtheorem{sledst}[theor]{\hspace{1cm}{\sc Corollary}}
\newtheorem{lemma}[theor]{\hspace{1cm}{\sc Lemma}}

\newtheorem{constr}[theor]{\hspace{1cm}{\sc Construction}}
\newtheorem*{utver*}{\hspace{1cm}{\sc Proposition}}
\theoremstyle{definition}

\newtheorem{defin}[theor]{\hspace{1cm}{\sc Definition}}
\newtheorem*{defin*}{\hspace{1cm}{\sc Definition}}
\newtheorem{exa}[theor]{\hspace{1cm}{\sc Example}}

\newtheorem{rem}[theor]{\hspace{1cm}{\sc Remark}}

\newcommand{\codim}{\mathop{\rm codim}\nolimits}

\newcommand{\Cr}{\mathop{\rm Cr}\nolimits}

\newcommand{\sing}{\mathop{\rm sing}\nolimits}
\newcommand{\sm}{\mathop{\rm sm}\nolimits}

\newcommand{\conv}{\mathop{\rm conv}\nolimits}

\newcommand{\Trop}{\mathop{\rm Trop}\nolimits}
\newcommand{\HTrop}{\mathop{\rm HTrop}\nolimits}

\newcommand{\dist}{\mathop{\rm dist}\nolimits}

\newcounter{idx}

\newcommand{\rotraise}[1]{
  \StrLen{#1}[\slen]
  \forloop[-1]{idx}{\slen}{\value{idx}>0}{
    \StrChar{#1}{\value{idx}}[\crtLetter]
    \IfSubStr{tlQWERTZUIOPLKJHGFDSAYXCVBNM}{\crtLetter}
      {\raisebox{\depth}{\rotatebox{180}{\crtLetter}}}
      {\raisebox{1ex}{\rotatebox{180}{\crtLetter}}}}
}

\renewcommand{\emph}[1]{{\it {\color{NavyBlue} #1}}}

\def\Hom{\text{Hom}}

\def\R{\mathbb R}

\def\Z{\mathbb Z}
\def\Q{\mathbb Q}
\def\C{\mathbb C}
\def\CC{({\mathbb C}^\star)}

\def\CP{\mathbb C\mathbb P}

\def\1{\mathbbm 1}
\def\sm{{\rm sm}\,}
\usepackage{xr-hyper}
\makeatletter

\newcommand*{\addFileDependency}[1]{%
\typeout{(#1)}%
\@addtofilelist{#1}
\IfFileExists{#1}{}{\typeout{No file #1.}}
}\makeatother

\newcommand*{\myexternaldocument}[1]{
\externaldocument[CR]{#1}
\addFileDependency{#1.tex}
\addFileDependency{#1.aux}
}

\myexternaldocument{crit}

\begin{document}

\begin{center}{\Large \sc Gauss--Bonnet for CSM classes of very affine manifolds, and tropical isotopies}

\vspace{3ex}

{\sc Alexander Esterov}
\end{center}

\begin{abstract}
The CSM class of a very affine manifold $U$ is represented by the rank drop locus of a general tuple of torus invariant 1-forms on it. This equality holds in the homology of any toric compactification $X\supset U$. It was proved for sch\"on $U$ by Huh, and later for all $U$ in the homology of $X=\CP^n$ by Maxim--Rodriguez--Wang--Wu, using Ginsburg's interpretation of CSM classes as Lagrangian cycles.

We deduce this identity in full generality from properties of affine characteristic classes,  
give an explicit sufficient condition of general position for the 1-forms, and use it to extend the identity to non-torus invariant 1-forms.
Along the way, we characterize very affine varieties of ML degree 0, and give a useful criterion for a family of varieties to have a constant tropical fan (``tropical isotopy theorem''). 

Our main idea is a technique to study very affine varieties without compactifying them. As a result we do not have to deal with singularities occuring at the boundary of the compactification, which would require less constructive methods (such as resolution of singularities). This is what enables us to give constructive genericity conditions in the Gauss--Bonnet theorem.
\end{abstract}

\tableofcontents
\section{Introduction}\label{sintro}

{\bf Huh's Gauss--Bonnet.} The Chern class of a compact complex manifold $U$, as a cycle in $H_{2k}(U)$, can be represented as the rank drop locus of general 1-forms $w_0,\ldots,w_{k}$ on $U$:  $$Z(w_0\wedge\ldots\wedge w_{k}):=\{x\in U\,|\,w_0,\ldots,w_{k}\mbox{ are linearly dependent at }x\},$$ provided that $\dim_\R Z(w_0\wedge\ldots\wedge w_{k})=2k$, and taking its components with natural multiplicities.
(We prefer to work with homology, but one can interchangeably use Chow rings.)

This Gauss--Bonnet equality fails for non-compact $U$. However, once such $U$ is {\it very affine}, i.e. admits a closed embedding into the algebraic torus $\CC^n$, the Gauss--Bonnet equality revives for generic torus-invariant 1-forms $w_i(x_1,\ldots,x_n)=c_1\frac{dx_1}{x_1}+\cdots+c_n\frac{dx_n}{x_n}$ and CSM classes of $U$, once $U$ is regarded as a semialgebraic subset of any toric compactification $X\supset\CC^n$.

\begin{theor}\label{th1}
Assume that the singular locus of a variety $U\subset\CC^n$ has dimension smaller than $r$, and $w_0,w_1,\ldots$ are generic invariant 1-forms on the smooth part. Then we have $(-1)^{\dim U-k} c^{SM}_k(U)=[\overline{Z(w_0\wedge\ldots\wedge w_{k})}]\in H_{2k}(X)$ for any $k\geq r$ and any toric $X\supset\CC^n$.

In particular, for $k=0$, the Euler characteristics $e(U)$ is the signed number of zeros of $w_0|_U$.
\end{theor}
Here $Z(\ldots)$ is defined as above for the smooth part of $U$, and $[\overline{Z(\ldots)}]$ is the fundamental cycle of its closure in $X$.
This was proved in \cite{huh} for $U$ sch\"on or $k=0$, and in \cite{maxim} for $X=\CP^n$. The methods of the first paper are specific for the sch\"on setting if $k>0$ (see Remark 3.10 therein) and involve Hironaka resolution of singularities of $\bar U\subset X$ for $k=0$. The second paper is based on an advanced technique of conic Lagrangian cycles.

We prove this theorem (in full generality) more constructively: while we do rely on transcendental methods at one point and thus stick to complex numbers, we avoid Hironaka resolutions or other constructions difficult to extend to arbitrary ground fields. Moreover, we actually never compactify $U$ and its ambient torus in the core arguments, except for one step in the proof of a disjoint auxiliarry statement (tropical isotopy theorem).

As a result, we give constructive genericity assumptions for the form $w:=w_0\wedge\ldots\wedge w_k$ to satisy Theorem \ref{th1}, which do not rely e.g. on the Hironaka or even Nash resolution of singularities of $\bar U$ in the toric boundary of $X$. This may facilitate extending the topic to arbitrary fields and arbitrary ambient reductive groups, as well as computing CSM classes of $U$ in terms of its ideal (using computer algebra systems or even by hands). 

\vspace{1ex}

{\bf The tropics.} The genericity conditions are given in terms of tropical fans of sets $Z(w)$. Recall that a $k$-dimensional very affine variety $V\subset\CC^n$ has the {\it tropicalization} $\Trop V$: it is a $k$-dimensional weighted polyhedral fan in the vector space $N:=M^*\otimes\R$ dual to the monomial lattice $M:=\Hom(\CC^n,\C^*)$, see e.g. \cite{ms} or \cite{mr} for an introduction. 

The correspondence $V\mapsto\Trop V$ is constructive: in an algebraic family of varieties $V_s,\, s\in S$, almost all varieties have the same tropical fan, equal to $\Trop V_\bullet$ over $\C(S)$, and referred to as the {\it generic} one. We denote $\Trop Z(w)$ for generic $w=w_0\wedge\ldots\wedge w_k$ by $\Trop_k U$. 

The weighted fans $\Trop_0 U
\subset\cdots\subset\Trop_{\dim U} U=\Trop U$ are called the {\it tropical characteristic classes} of $U$ and are the main characters in our story. We set $\Trop_j U=0$ for $j>\dim U$.

\begin{rem}
In the introduction, we look at the fans $\Trop_k U$ 
only for $k>\dim\sing U$. However, in Section \ref{saff}, we relate them to another entity, the {\it affine characteristic class} of $U$, which actually generalizes $\Trop_k U$ to every $k$ (for arbitrarily singular $U$).
\end{rem}

The fan $\Trop Z(w)$ represents its generic value $\Trop_k U$ if the form $w$ satisfies the same genericity condition as the one promised for Theorem \ref{th1}. This will allow to express the fan $\Trop_k U$ (and thus the genericity condition itself) constructively in terms of the ideal of $U$. 

\vspace{1ex}

{\bf The genericity conditions}  
will be given in terms of the mutual position of the plane $\ker w\subset T_1\CC^n$ and the fans $\Trop_j U\subset N$ (note that the two ambient spaces are naturally isomorphic). 
Let $p$ be the projection $N\to N/\ker w$. 
\begin{defin}
1. The projection $p$ (or $\ker w$) is said to {\it immerse} a 
polyhedral complex $K\subset N$ of pure dimension, if its restriction to a suitable open dense subset of $K$ is injective.

2. The projection $p$ is said to {\it submerse} a  
polyhedral complex $K\subset N$ of pure dimension, if all fibers of $p|_K$ have dimension not higher than the generic one (i.e. $\dim K-\dim$(range of $p)$).
\end{defin}
\begin{exa}
The complexes of dimension 0 and 1 on the picture are not submersed or immersed (respectively) by the projection $p$, but we can fix it by removing the blue elements. 
\begin{center}
\includegraphics[scale=0.6]{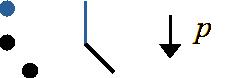}
\end{center}
\end{exa}
\begin{theor}\label{tropth} $\,$ 

1. In the setting of Theorem \ref{th1}, we have $\dim Z(w)=k$ and $p\Trop Z(w)=p\Trop_k U$ if  $$
\Trop_qU\mbox{ is submersed by }p\mbox{ for }q>k,\eqno{(*)}$$ and moreover $\Trop Z(w)=\Trop_k U$ if additionally $$\Trop_k U\mbox{ is immersed by }p.\eqno{(**)}$$ 

2. Theorem \ref{th1} holds for any form $w$ satisfying $(*)$ and $(**)$. 
\end{theor}
\begin{rem} \label{algorithm} 1. Both $(*)$ and $(**)$ are conditions of general position on the form $w$. Once they are satisfied, we see that the fan $\Trop Z(w)$ represents the tropical characteristic class $\Trop_k U$, and the class $[Z(w)]$ represents the CSM class $c^{SM}_k(U)$.

2. In particular, $(-1)^{\dim U}e(U)$ is the weight of the fan $\Trop_0 U$ and the weighted number of zeroes of any 1-form $w|_U$ satisfying condition (*) (because (**) is trivial for $k=0$).
\end{rem}

{\bf Constructivity.} The theorem allows to reconstruct the classes $\Trop_k U$ inductively for all $k$ from $n$ to 0. This is not tautological, because there is an apparent circular dependence: the condition $(**)$ for $\Trop Z(w)$ to represent the sought fan $\Trop_k U$ depends on $\Trop_k U$ itself. 

To keep things simple, we show how to overcome this apparent circular dependence in the first non-trivial case, when $U\subset\CC^3$ is a surface.

\begin{constr}\label{constr1}

\begin{adjustwidth}{4ex}{0ex}
(i) The classes $\Trop_3 U=0$ and $\Trop_2 U=\Trop U$ are known.

\noindent (ii) To find $\Trop_1 U$, choose a 2-form $w$ so that the line $\ker w$ intersects the known 2-dimensional fan $\Trop_2 U$ by $0$. Then $\Trop_1 U$ belongs to the 1-dimensional fan $K:=(\Trop_2 U)\cap p^{-1}p\Trop Z(w)$: here we use 
the first equality of the theorem. 

\noindent (iii) Choose another 2-form $\tilde w$ so that the projection along its kernel immerses $K$ (and thus $\Trop_1 U$, which is still unknown at this point, but is contained in $K$), and $\ker\tilde w\cap \Trop_2 U=\{0\}$. Then $\Trop_1 U=\Trop Z(\tilde w)$: here we use the second equality of the theorem.

\noindent (iv) Choose a 1-form $w$ so that the plane $\ker w$ intersects $\Trop_2 U$ by at most a 1-dimensional set and $\Trop_1 U$ by 0, then $\Trop_0 U=\Trop Z(w)$.
\end{adjustwidth}

\end{constr}

This algorithm assumes that we know the fan $\Trop Z(w)$ for a rational form $w$ and a variety $U$ given by its ideal $I$. This is so: first, $Z(w)$ is the critical locus of the map $p|_U$, so its ideal is constructive from $I$. Second, the tropical fan of a variety is constructive from its ideal.

Once we have a machine computing the tropical fan of a critical locus of a projection for a variety given by its ideal, the results of our paper allow to both compute the tropicalized rank drop locus $\Trop Z(w)$ for a generic form $w$ and to tell it from a non-generic one.

\vspace{1ex}

{\bf Weights.} In this context, $Z(w)$ should be regarded as an algebraic cycle, whose components have natural positive weights. The weight at a generic point $z\in Z(w)$ is defined as follows: 

-- the number of branches of the (reduced) set $Z(\tilde w)$ near $z$ for a geneirc form $\tilde w$ close to $w$;

-- the Milnor number of an ICIS $m_0=\cdots=m_k=0$ on $U$ such that $w=dm_0\wedge\ldots\wedge dm_k$.

See Proposition \ref{ppencil}.1 for the equivalence of the two definitions. Note that the second one admits an algebraic expression (L\^e--Greuel formula) in terms of the local ideal of $U$ at $z$.

These weights are not important for Theorem \ref{th1}, as they equal 1 for sufficiently generic $w$. However they may fail to be 1 under particular genericity conditions (*,**) of 
Theorem \ref{tropth}.

Thanks to our assumption $k>\dim\sing U$, we can interchangeably discuss $k$-cycles on $U$ and on its smooth part $\sm U$. In particular, this applies to the $k$-cycle $Z(w)$.

\vspace{1ex}

We conclude the introduction with some applications of our results, and some key ideas of the proof (which may be of independent interest).

\vspace{1ex}

{\bf Computing the CSM class $c^{SM}_k(U)$.} The Fulton--Sturmfels isomorphism \cite{fs} is a combinaotrial model for the homology ring of a smooth toric compactification $X_\Sigma\supset\CC^n$ defined by a toric fan $\Sigma$. Namely, it is the graded ring isomorphism
$${\rm HTrop}_\Sigma:H_k(X_\Sigma)\to H_k(\Sigma):=\{\mbox{tropical fans supported on the }k\mbox{-skeleton of }\Sigma\},$$
uniquely defined by compatibility with the tropicalization: if a $k$-dimensional variety $V\subset X_\Sigma$ avoids orbits of codimension higher than $k$, then
${\rm HTrop}_\Sigma$ sends the fundamental class $[V]$ to the tropical fan $\Trop V\cap\CC^n$. If $F\subset H_k(\Sigma)$, we say that the fan $\Sigma$ 
{\it refines} $F$.

Applying ${\rm HTrop}^{-1}$ to the equality in Theorem \ref{th1}, we get the following from Theorem \ref{tropth}.2.
\begin{theor}\label{tropcsm}
Choose a smooth toric compactification $X\supset\CC^n$ whose fan refines $\Trop_k U$. Then the CSM class $c^{SM}_k(U)\in H_k(X)$ equals  $(-1)^{\dim U-k}{\rm HTrop}^{-1}\Trop_k U$, i.e. is represented by the fan $(-1)^{\dim U-k}\Trop_k U$ in the Fulton--Sturmfels model of $H_k(X)$.

In particular, the Euler characteristics $e(U)$ is the weight of the 0-dimensional fan $\Trop_0 U$.
\end{theor}
Thus our results on $\Trop_k U$ apply to $c^{SM}_k(U)$.

\vspace{1ex}

{\bf Non-proper non-invariant Poincare--Hopf.}
Checking the genericity conditions (*) and (**) for a seemingly far-from-generic 1-form $dy/y$ on the graph $\{g(x)=y\}$ of a polynomial $g:U\to\C$, we get a version of Theorem \ref{th1} for generic non-invariant 1-forms. 
\begin{theor}\label{thnoninv}
For a finite set of monomials $A\subset M$, let every face $0\in\Gamma\subset\conv\{0\}\cup A$ have no exterior normal vectors in the fan $\Trop_{dim\Gamma} U$. Then a generic linear combination $g$ of monomials from $A$, regarded as a polynomial function on $U$, has a finite number of critical points, 
equal to $(-1)^{\dim U}e(U\setminus\{g=0\})$ and to $\sum_{k\geq 0}(\Trop A)^k\Trop_k U$.
\end{theor}
Here $\Trop A$ is the tropical fan of the polytope $\conv A$ (i.e. the corner locus of its support function), and a product of tropical fans is evaluated as their tropical intersection number.

Note that most of sets $A$ suit this theorem for given $U$ (e.g. all sets of size $a$ outside a hypersurface in $M^a$); c.f. another occurrence of such ``generic Newton polytopes'' in \cite{khgelf}.

The case when $A$ consists of one monomial $m$ is exactly Theorems \ref{th1} and \ref{tropcsm} for 1-forms.

\vspace{1ex}

{\bf Tropical isotopy theorem} claims that a tropicl fan is constant in a family, if its sufficiently ``ample'' projection is. It is an algebraic geometry implication of a tropical fan being immersed.

\begin{theor}[see Section \ref{sisot}]\label{tropisot}
Let $V\subset\CC^n\times\C$ be of pure dimension $k+1$, and its fibers $V_t,\, t\in\C^*$, under the coordinate projection $\CC^n\times\C\to\C$ of pure dimension $k$. Let $\Trop V_\bullet\subset\R^n$ be the tropical fan of the generic fiber, and $p:\CC^n\to\CC^r$ an epimorphism such that $p_*\Trop V_\bullet=p_*\Trop V_0$. If $\Trop V_\bullet$ is immersed by $p$, then $\Trop V_\bullet=\Trop V_0$.
\end{theor}
In particular, the fan $\Trop V_t$ is the same for all small $t$, hence the name "tropical isotopy".

\vspace{1ex}

{\bf Varieties of ML degree zero.} The proof of the main results relies on a characterization of such varieties. Thoughit is simpler than for ML degree one in \cite{huh2}, it is surprisingly difficult to find it in the literature. So we give it here with a proof. 
\begin{theor}[see Section \ref{sgauss}]\label{thclassif}
Consider a (not necessarily smooth) variety $V\subset\CC^n$.

1. If a generic torus invariant 1-form has no zeroes when restricted to the smooth part $\sm V$, then $V$ is the preimage of a variety under a homomorphism $\CC^n\to\CC^{n-1}$. 

2. More generally, if $\Trop_k\sm V=0$ for $k<r$, then $V$ is a preimage under $\CC^n\to\CC^{n-r}$.
\end{theor}
\begin{sledst}
A manifold in $\CC^n$ has Euler characteristics 0 if and only if it preserved by a proper subtorus of $\CC^n$.
\end{sledst}
Interestingly, sch\"on varieties satisfying this corollary are exactly those that are not log canonical \cite[Theorem 3.1]{fromhubsch}. It would be interesting to extend this to all very affine manifolds.

\vspace{1ex}

{\bf The structure of the paper.} In Section \ref{saff}, we relate the topic to affine characteristic classes of very affine varieties (with arbitrary singularities) constructed in \cite{jems}. From this relation (Theorem \ref{thropaff}) and Theorem \ref{tropth}.1, the other results of the introduction follow immediately (see the end of Section \ref{saff}). Besides that, we use this relation to constructively compute affine characteristic classes (Corollary \ref{sledaff}) and to express the classes $\Trop_\bullet U_1\cap U_2$ in terms of $\Trop_\bullet U_i$, if $U_1$ and $U_2$ are in general position (Corollary \ref{sledtrop}).

Sections \ref{sdiscr} and \ref{sgauss} are devoted to the proof of the key Theorems \ref{thropaff} and \ref{tropth}.1. It is convenient to prove them simultaneously by induction on the increasing $n$ (the dimension of the ambient torus) and then decreasing $k$ (the dimension of the characteristic class). However, even if we start from a smooth variety, this induction reduced the statement to smaller dimensional varieties with at least isolated singularities. So, in Section \ref{sdiscr}, we relax the assumptions of Theorems \ref{thropaff} and \ref{tropth}.1 to admit varieties with isolated singularities even for $k=0$, and then finish the inductive argument in this higher generality in Section \ref{sgauss}. 

A key tool in this argument is the tropical isotopy theorem, which we prove in Section \ref{sisot}. Finally, Section \ref{sdiscus} contains some discussion and open questions.

\vspace{1ex}

The author is grateful to Lauren\c tiu Maxim for enlightening discussions.

\section{Affine characteristic classes}\label{saff}

{\bf Affine characteristic class.} Tropical fans in $\R^n$ form a graded ring $K=K_0\oplus\ldots\oplus K_n$, in which $K_j$ is the vector space of codimension $j$ tropical fans, and multiplication is the tropical intersection product (see e.g. the original paper \cite{fs} or \cite{mr}). 
By \cite{jems}, there is is a unique map 
$$\mbox{subvarieties of }\CC^n\to K,\; U\mapsto \langle U\rangle=\langle U\rangle_0+\cdots+\langle U\rangle_n,$$
with the properties: (i) $\langle U\cdot V\rangle=\langle U\rangle\cdot\langle V\rangle$, i.e. $\langle U\cdot V\rangle_k=\sum_j\langle U\rangle_j\cdot\langle V\rangle_{k-j}$ for $U,V\subset\CC^n$; 

\vspace{1ex}

(ii) $\langle U\times V\rangle=\langle U\rangle\times\langle V\rangle$, i.e. $\langle U\times V\rangle_k=\sum_j\langle U\rangle_j\times\langle V\rangle_{k-j}$ for $U\subset\CC^n$ and $V\subset\CC^m$; 

\vspace{1ex}

(iii) $\langle U\cup V\rangle+\langle U\cap V\rangle=\langle U\rangle+\langle V\rangle$. Thanks to this, given a {\it constructive function} $U:\CC^n\to\Z$ (i.e. a linear combination $\sum_i a_i\1_{U_i}$ of the indicator functions of varieties $U_i$), we can correctly define its characteristic class $\langle U\rangle:=\sum_i a_i\langle U_i\rangle$. This allows us to formulate the next property. 

\vspace{1ex}

(iv) For a projection of tori $p:\CC^n\to\CC^m$, we have $p_*\langle U\rangle=\langle p_*U\rangle$, where the {\it MacPherson image} $p_*$ of a variety $U\subset\CC^n$ is the constructive function $\CC^m\to\Z,\, y\mapsto e\Bigl(U\cap p^{-1}(y)\Bigr)$.

\vspace{1ex}

(v) For a smooth toric compactification $X\supset\CC^n$ whose toric fan refines the fan $\langle U\rangle_k$, the Fulton--Sturmfels isomorphism ${\rm HTrop}$ identifies this fan with the CSM class $c^{SM}_{n-k}(U)\in H_\bullet(X)$. 

\vspace{1ex}

Note that we prefer to index $\langle U\rangle_k$ with its codimension (in contrast to indexing $\Trop_k U$ and $c^{SM}_k U$ with the dimension), in order to simplify summation indices in (i-ii).

\begin{exa}
1. For $U\subset\CC^n$, the number $\langle U\rangle_n$ is the Euler characteristics $e(U)$.

2. If $\codim U=k$, then $\langle U\rangle_i=0$ for $i<k$, and $\langle U\rangle_k$ is the tropicalization $\Trop U$.
\end{exa}

\begin{rem}
Property (v) was formulated in \cite{jems} for a ``sufficiently fine'' toric compactification $X$, not specifying what particular toric fans are sufficiently fine for it. Let us prove that the same holds for any smooth toric variety $X'$, whose fan refines $\langle U\rangle_k$.

For this, choose a smooth toric variety $Y$ refining both of them: $p:Y\to X$ and $p':Y\to X'$. We know that $\HTrop p^*\alpha=\HTrop\alpha$ for any $\alpha\in H_\bullet(X)$ and $\HTrop p'_*\beta=\HTrop\beta$ for any $\beta\in H_\bullet(Y)$ under an obviously necessary assumption that $\HTrop\beta\in\HTrop H(X')$ (i.e. the toric fan of $X'$ refines $\HTrop\beta$). Denote $U$ as a subset of $X, X', Y$ by $U_X, U_{X'}, U_Y$, and onsider 

-- the first equality for $\alpha=c^{SM}(U_X)\in H_\bullet(X)$; 

-- the second equality for $\beta=c^{SM}(U_Y)\in H_\bullet(Y)$;

-- the equality $p^* c^{SM}(U_X)=c^{SM}(p^{-1}U_X)=c^{SM}(U_Y)\in H_\bullet(Y)$; 

-- the equality $p'_*c^{SM}(U_Y)=c^{SM}(p'_*U_Y)=c^{SM}(U_{X'})\in H_\bullet(X')$.

These four equalities combine to the sought one.
\end{rem}

{\bf The relation between affine and tropical characteristic classes.}
\begin{theor}\label{thropaff}
For every $k<n-\dim\sing U$, we have $\langle U\rangle_{k}=(-1)^{\dim U+k-n}\Trop_{n-k} U$.
\end{theor}
Unlike Theorem \ref{th1} or property (v) of affine classes, this is an equality of tropical fans, independent of any choice of a troic compactification. We apply this equality in both directions: 

-- to facilitate constructing $\langle U\rangle_{k}$ from the ideal of $U$;

-- to compute $\Trop_k U\cap g\cdot V$ for given varieties $U$ and $V\subset\CC^n$ and a generic shift $g\in\CC^n$.

\vspace{1ex}

The class $U\mapsto \langle U\rangle$ is explicitly constructed in \cite{jems} from the tropical fan of the conormal cycle of $U$. However, this construction is of little help for practical computations, because passing to conormal cycles greatly increases the dimension of the variety and its ambient torus.

\begin{sledst}\label{sledaff}
Given a set $U\subset\CC^n$ of pure dimension $m$, choose a polynomial $f$ vanishing on $\sing U$, but not on components of $U$. Then, for the smooth set
$$U':=\{(x,t),|\,x\in U,\,t\cdot f(x)=1\}\subset\CC^n\times\C^*\xrightarrow{\pi}\CC^n$$
and the $(m-1)$-dimensional set $S:=U\cap\{f=0\}$, we have:
$$\langle U\rangle_k=\langle S\rangle_k+(-1)^kp_*\Trop_{n-k} U'.$$
\end{sledst}
This is a way to construct $\langle U\rangle_k$, because both terms on the right are constructive: the first one by induction on $m$, and the second one as explained in the introduction. 

\begin{sledst}\label{sledtrop}
For varieties $U$ and $V\subset\CC^n$ and a generic shift $g\in\CC^n$, we have $$\Trop_k U\cap gV=\sum_j \Trop_j U\cdot \Trop_{n-k-j} V.$$
In particular, for $k=0$ we get the Euler characteristics: $$e(U\cap gV)=\sum_j \Trop_j U\cdot \Trop_{n-j} V.$$
\end{sledst}
This is not trivial to see by intersecting particular rank drop loci $Z(\ldots)$ on $U$ and $V$.

\vspace{1ex}

{\bf Some reductions.} In the presented results, the central ones are Theorems \ref{tropth}.1 and \ref{thropaff}. From them, we now deduce the rest (Theorems \ref{tropcsm}, \ref{tropth}.2, and thus \ref{th1}).  
Note that in the introduction we offered Theorem \ref{tropcsm} as a corollary of Theorem \ref{th1}. We now notice that it is rather a corollary of Theorem \ref{thropaff} and property (v) of affine characteristic classes. This allows to apply it the other way around, deducing Theorem \ref{tropth}.2 (and thus \ref{th1}).

\vspace{1ex}

{\it Deducing Theorem \ref{tropth}.2 from \ref{tropcsm} and \ref{tropth}.1.} If the toric fan of the variety $X$ refines $\Trop_k U$, then both $[Z(w)]$ and $c^{SM}_k(U)$ are sent to $\Trop_kU$ by Fulton--Sturmfels isomorphism $\HTrop$: 

-- the first one by Theorem \ref{tropth}.1 and the identity $\HTrop [V]=\Trop V$ for any variety $V\subset\CC^n$ whose fan $\Trop V$ belongs to $\HTrop H(X)$ (i.e. is refined by the fan of $X$);

-- the second one by Theorem \ref{tropcsm}.

From this, we see $[Z(w)]=c^{SM}_k(U)=\HTrop^{-1}\Trop U$ once the fan of $X$ refines $\Trop_k U$.

\vspace{1ex}

Now consider an arbitrary smooth toric variety $X$, choose a smooth toric blow up $p:Y\to X$ such that the toric fan of $Y$ refines $\Trop_k U$. The map $p_*:H(Y)\to H(X)$ sends $[Z(w)]$ to $[Z(w)]$ (as $p$ sends $Z(w)$ to $Z(w)$) and $c^{SM}_k(U)$ to $c^{SM}_k(U)$ (as $p_*c^{SM}_k(U)=c^{SM}_k(p_*U)$, and $p_*U=U$). Thus the known equality $[Z(w)]=c^{SM}_k(U)$ in $H(Y)$ implies the same in $H(X)$ (though $\HTrop$ may not send these equal elements of $H(X)$ to the fan $\Trop_k U$ anymore).
$\hfill\square$

\section{Discriminants of projections}\label{sdiscr}

{\bf The Bertini discriminant.} Denote the additive topological Euler characteristics by $e$.

\begin{defin}\label{defcr}
1. A function $\varphi$ on an analytic space $U$ is said to be {\it constructive}, if it is a linear combination of indicator functions $\1_{U_i}$ of finitely many subvarieties $U_i\subset U$.

2. If $\varphi$ can be represented as $d+\sum_{i:\codim U_i\geq k} d_i\cdot \1_{U_i}$, then its (signed) {\it codimension $k$ discriminant} $D^{k}_\varphi$ is defined as the codimension $k$ cycle $(-1)^k\sum_{i:\codim U_i=k} d_i\cdot U_i$.

3. For a regular map $p:U\to V$, its {\it Milnor function} $\mu_p$ is defined at every point $x\in U$ as the Euler characteristics of its Milnor fiber, i.e. $e(W\cap\{p=v\})$ for a small metric neighborhood $W$ of $x$ and generic $v\in V$ close to $p(x)$ (relatively to the size of $W$).

4. If the constructive function $\mu_p$ has the discriminant of codimension $\dim U-\dim V$, then it is called the {\it critical locus} $\Cr p$, and the divisor $p_* \Cr p$ in $V$ the {\it Bertini discriminant} $B_p$.
\end{defin}

\begin{rem}\label{rembert}
1. We speak of analytic spaces rather than algebraic varieties, because at times we apply these notions to the preimage of a variety under $exp:\C^n\to\CC^n$. In particular, for any invariant form $w=w_0\wedge\ldots\wedge w_k$ in the\ torus (even if $\ker w$ is an irrational subspace in $T_1\CC^n$), we can define $\Cr w|_U$ as $\exp(\Cr p|_{\ln U})$, where $p$ is the projection $\C^n\to\C^n/\ker w$.

2. The definition of the Milnor function encrypts a standard statement: there exists a Zariski open set $\tilde V\subset V$ and an integer $\mu(x)$ constructively depending on $x\in U$, such that, for any smooth metrics on $U$ and $V$ we have: $\forall x\in U\; \exists \epsilon_0>0\; \forall \epsilon<\epsilon_0\; \exists \delta>0\; \forall y\in\tilde V\,\,:\dist(y,p(x))<\delta\,\Rightarrow\,(-1)^{\dim U-\dim V}\bigl( e(W_\epsilon\cap\{p=y\})-1\bigr)=\mu(x)$, where $W_\epsilon$ is the metric $\epsilon$-neighborhood of $x$.
The words ``generic'', ``sufficiently small'' and ``close'' in the definition should be understood in the sense of this statement.

3. If $U$ and $V$ are smooth, $\mu_p(x)$ is the usual Milnor number for an isolated singularity of a fiber of $p$, which is a complete intersection. In particular, the Milnor number is positive and can be computed algebraically by the L\^e-Greuel formula \cite{le}, \cite{greuel}.

4. In particular, the divisor $B_p$ is effective once $\dim\sing U<\dim V-1$.

5. Otherwise, it may be not effective: for instance, a component $C$ of the set $\sing U$ may have the same dimension as $\Cr p$ but have weight 0 in this algebraic cycle. Nevertheless, it is still clear how to compute this weight algebraically: for a generic map $\tilde p$ near a generic point $x\in C$, the value $\mu_{\tilde p}(x)$ expresses in terms of the Euler obstruction of $U$ at $C$ (given by the L\^e--Teissier algebraic formula \cite{lt}), and the difference 

$\mu_{p}(x)-\mu_{\tilde p}(x)$ for arbitrary $p$ is expressed algebraically by the Ebeling--Gusein-Zade--Seade homological index \cite{egs}.

6.  
In what follows, when we apply set-theoretic operations like $\subset$ or $\cap$ to $\Cr p$ or $B_p$, we understand them as sets $\Cr p:=\sing U\cup\{$critical points of $p$ on $\sm U\}$ and $p(\Cr p)$.  

We should be cautious here, as they may differ from the support sets of the respective algebraic cycles, see (5) above. However, this subtlety disappears once we assume $\dim\sing U<\dim V-1\ne 0$ (as we did in the introduction).

7. Similarly $\Trop \Cr p$ is understood as a polyhedral complex or a weighted fan depending on whether we apply set theoretic or intersection theoretic operations to it: some cones of the complex $\Trop \Cr p$ may have zero weights in the tropical fan $\Trop \Cr p$ if $\dim\sing U=\dim V-1$.

8. While the critical cycle is defined for any $p$, what follows has little sense unless $\dim\sing U\leq\dim V$, and $p$ restricted to the $k$-dimensional part of $\sing U$ has finite degree. So we always assume the first of these two conditions whenever we discuss $B_p$, and both conditions for $\Cr p$. This is important for instance for the statmenet of the following lemma.
\end{rem}

{\bf Fubini type formulas for Bertini discriminants.}

\begin{lemma}\label{lfub}
Let $(W,\alpha)\xrightarrow{p}(\C^2,0)\to(\C,0)$ be germs of maps such that their composition $t$ has reduced fiber $\{t=0\}$, and the cycle $\Cr t$ is defined near $\alpha$ (thus supported at $\alpha$). Then 
$$\Cr {p}|_{t=0} -\{t=0\}\cdot\Cr p+\Cr t =0,$$
where the multiplication is the intersection product. This equality implies that that the first term is defined (i.e. $dt$ does not vanish outside $\alpha$, and $dp/dt$ does not vanish on $\{t=0\}\setminus\{\alpha\}$) if and only if the second term is (i.e. $\Cr p $ is defined and intersects $t=0$ by $\alpha$).
\end{lemma}

\begin{proof} Let us choose a sufficiently small metric neighborhood $U$ of $\alpha\in W$, a sufficiently smaller square neighborhood $V$ of $0\in\C^2$, and a number $t_0\ne 0$ of sufficiently smaller absolute value. Then we have the following relations.

1. $e(U\cap\{t=t_0\})=e(U\cap\{t=t_0\}\cap p^{-1}V)+e(U\cap\{t=t_0\}\setminus p^{-1}V)$ by the additivity of the Euler characteristics. We now study each of the two terms on the right.

2. $e(U\cap\{t=t_0\}\cap p^{-1}V)=e($generic fiber of $p)-(-1)^{\dim W}\{t=t_0\}\cdot\Cr p$, counting the left hand side fiberwise over $p$. Up to the sign $(-1)^{\dim W}$, the right hand side can further be simplified to $\Cr 
 {p}|_{t=0}+1-\{t=0\}\cdot\Cr p $, because $p$ has the same generic fiber as $p|_{t=0}$.

3. $e(U\cap\{t=t_0\}\setminus p^{-1}V)=e(U\cap\{t=0\}\setminus p^{-1}V)=e(U\cap\{t=0\})-e(U\cap\{t=0\}\cap p^{-1}V)=1-1=0$, because $U\cap\{t=t_0\}\setminus p^{-1}V\cong U\cap\{t=0\}\setminus p^{-1}V$.

The three equalities sum up to the statement of the lemma.
\end{proof}

\begin{utver}\label{pfub}
Let $(W,\alpha)\xrightarrow{p}(V\times\C,0)\to(\C,0)$ be germs of maps such that their composition $t$ has reduced fiber $\{t=0\}$, and its Milnor function $\mu_t$ has discriminant $D$ of codimension $\dim W-\dim V$ in $\{t=0\}$ near $\alpha$. Then the cycle $$\Cr {p}|_{t=0} -\{t=0\}\cdot\Cr p+D$$
is vertical, i.e. the image of its underlying set in $V$ has codimension higher than 1. This equality implies that that the first term is defined if and only if the second term is.
\end{utver}
This is from the preceding lemma, by restricting to a small nrighborhood of a point $x\in \{t=0\}$ near $\alpha$ in $W$ and to a generic germ of a smooth curve at $p(x)$ in $V$. Importantly, $D$ does not depend on $p$ here, only on its composition $t$. In particular, if $t$ is locally trivial, then $D=0$, and we get the following conservation law for critical loci.
\begin{sledst}\label{conserv}
For a family of maps $p_t:U\to V$, if the critical locus $\Cr p_0$ is defined, then the family of cycles $\Cr p_t$ is defined for small $t$, and the cycle $\lim_{t\to 0}\Cr p_t$ (in the sense of the definition below) differs from $\Cr p_0$ by a vertical cycle.
\end{sledst}
\begin{defin}\label{defconserv}
A family of cycles $C_t$ in an analytic space $U$ is {\it numerically continuous} at $t=0$, if, for any curve germ $\gamma\subset U$ properly intersecting $C_0$ we have $C_0\cdot\gamma=C_t\cdot\gamma$ for small $t$. \newline We shall take the liberty to express this mnemonically as $C_0=\lim_{t\to 0} C_t$.
\end{defin}

{\bf Main results for mildly singular varieties.} We restate the main results from the introduction for a $(k+1)$-form $w$ on $U\subset\CC^n$ in a slightly more general case: we allow isolated singularities of $U$ even for $k=0$. This is because we need it in the proof: even to deduce the results for a smooth $U$, we will have to apply them inductively to varieties of smaller dimension with isolated singularities. 

At the same time, we reduce the initial generality by assuming that the form $w$ has rational coefficients. This implies that it is the pull back of the invariant volume form under some torus projection $p:\CC^n\mapsto\CC^{k+1}$: in coordinates, $p=(p_0,\ldots,p_k)$ and $w=\frac{dp_0}{p_0}\cdots\frac{dp_k}{p_k}$. In this setting, we can interpret the rank drop locus $Z(w)$ as the critical locus $\Cr p$.
\begin{theor}\label{thfull}
Consider a projection $p:U\subset\CC^n\mapsto\CC^{k+1}$ with $\dim\sing U<k$ if $k>0$ or with at most isolated singulaities for $k=0$. 

1. We have $\langle U\rangle_j=(-1)^{\dim U+j-n}\Trop_{n-j} U$ for $j\leq n-k$.

2. If the projection $dp$ submerses the fan $\Trop_j U$ for $j> k$, then $dp_*\Trop \Cr p=dp_*\Trop_k U$.

3. If in addition $k=0$, then, in particular, the Euler characteristics $e(U)$ is the sum of the Milnor numbers $\mu_p(x)$ of $p$ at its critical points $x$ (including the singular points of $U$).

4. If $k$ is arbitrary, but in addition $dp$ immerses $\Trop_k U$, then $\Trop \Cr p=\Trop_k U$.
\end{theor}

Note that (3) is immediate from (1) and (2). The proof of (1,2,4) occupies the rest of this section and the next one. Before proceeding to it, we show that Theorem \ref{thfull} implies the rest of our results (all of which are reduced to Theorems \ref{thropaff} and \ref{tropth}.1 in Section 2\ref{saff}).

\vspace{1ex}

{\it Deducing Theorems \ref{thropaff} and \ref{tropth}.1 from \ref{thfull}.} Theorem \ref{thfull}.1 coincide with the sought Theorem \ref{thropaff}. Theorem \ref{thfull}.2\&4 coincides with the sought Theorem \ref{tropth}.1 for forms $w$ defined over $\Q$. So it remains to extend it to irrational forms. Such forms are parameterized by the points of the Grassmannian $G:=GR(n,n-k-1)$ by sending $w$ to the plane $\ker w\subset T_1\CC^n$. Consider the closure $S\subset G$ of all forms $w$ with non-generic values of the fan $\Trop Z(w)\ne\Trop_k U$ or its image $dp_*\Trop Z(w)\ne dp_*\Trop_k U$. Resolving singularities of the closure of $U$ in $\CP^n\supset\CC^n$, one can show that $S$ is a finite union of Schubert varieties defined over $\Q$. (For $k=0$, this is proved in \cite{huh}, and extends to $k>0$.) We need to show that this exceptional union of Schubert cells is covered by the Schubert varieties described by failing to satisfy conditions (*) and (**) of Theorem \ref{tropth}.1. Towards the contradiction, assume that some Schubert variety $C\subset S$ does not entirely consist of forms failing to satisfy (*) and (**). Since it is defined over $\Q$, we can take a generic form $w\in C$ with rational coefficients 
and notice that it contradicts Theorem \ref{thfull}.4 or 2 respectively: $w$ being generic in $S$ means failing to satisfy the conclusion of Theorem \ref{thfull}.4 or 2, while being generic in $C$ is equivalent to satisfying its assumption. $\hfill\square$

\vspace{1ex}

This argument suggests that extending our results to arbitrary fields will be significantly simpler for rational forms than for arbitrary ones.
Our proof of Theorem \ref{thfull} is based on comparing the Bertini discriminant $B_p$ to another divisor which we introduce now.

\vspace{1ex}

{\bf The Euler discriminant.} If the map $p:U\to V$ is proper with $\dim\sing U<\dim V$, the point $v$ is generic in a weight $\mu$ component of the Bertini discriminant $B_p$, and $\tilde v$ is generic in $V$, then we have $e p^{-1}(\tilde v)-e p^{-1}(v)=\mu$. We slightly strengthen this statement as follows.
\begin{defin} 1. The {\it McPherson pushforward} of $W\subset U$ under a morphism of algebraic varieties $p:U\to V$ is the function $(p_* W):V\to\Z$, sending $x\in V$ to $e(W\cap p^{-1}x)$.

2. The codimension 1 discriminant of the constructive function $p_*U$ (Definition \ref{defcr}.2) is called the {\it Euler discriminant} of $p$, denoted by $E_p$.
\end{defin}
\begin{utver}
For proper $p:U\to V$, the Euler discriminant equals the Bertini discriminant once the latter is defined: $E_p=B_p$.
\end{utver}
If $p$ is not proper, such as our composition $U\hookrightarrow\CC^n\mapsto\CC^q$, then $E_p$ still has functorial properties similar to Proposition \ref{pfub} for Bertini discriminants (see e.g. \cite{adv}), but the difference $E_p-B_p$ is non-trivial due to contributions from the infinity. These contributions are computed in \cite{adv} in case $U$ is a generic complete intersection. However, in the general setting they are complicated: see the literature surveyed in \cite{adv}, and e.g. \cite{ed24} and \cite{bd24} for recent advances concerning Euler and Bertini discriminants respectively. All the more fortunate looks the following observation.
\begin{theor}\label{thproj0}
1. Given a variety $U\subset\CC^n$ with $\dim\sing U\leq k$, for generic $p:\CC^n\mapsto\CC^{k+1}$, the Bertini discriminant of $p|_U$ is defined and equals the Euler discriminant.

2. A sufficient genericity assumption is that $dp$ submerses the tropical fans $\Trop_j U,\, j>k$. 
\end{theor}

{\it The proof assuming Theorem \ref{thfull}.3 for smaller $n$ and \ref{thfull}.1 for larger $k$.} 
I. We take $p$ satisfying the condition of (2), and prove (1) outside the non-properness locus $S$ of $p|_{\sing U}$. More specifically, we prove $B_p=E_p$ for $p$ over $\CC^{k+1}\setminus S$. Since $\codim S>1$ by the condition $\dim\sing U\leq k$, the equality $B_p=E_p$ outside $S$ implies the same on the whole torus $\CC^{k+1}$.

II. By the combinatorial Lemma \ref{lindtransv} below, appending $p:\CC^n\to\CC^{k+1}$ to $q:=(p,m):\CC^n\to\CC^{k+2}$ with a generic character $m:\CC^n\to\C^*$, the map $dq$ submerses $\Trop_i U$ for $i>k+1$ and immerses $\Trop_{k+1} U$. 
Thus Theorem \ref{thfull}.1 applies and gives $\Trop\Cr q=\Trop_{k+1} U$. Since the latter is assumed to be submersed by $dp$, the projection of $\Trop\Cr q$ along $\ker dp$ is proper, all the more so the restriction $p|_{\Cr q}$ over $\CC^{k+2}\setminus($non-properness locus of $q|_{\sing U})$, and all the more $p|_{\Cr p}$ over $\CC^{k+1}\setminus S$.

Note that we have to exclude non-properness loci here, because some components of $\sing U$ may be not contained in the support sets of the algebraic cycles $\Cr p$ or $\Cr q$, thus their tropical fans may be not contained in the tropical fans of these cycles.

\vspace{1ex}

III. An embedding $C:(\C,0)\to(\CC^{k+1}\setminus S,x_0)$ defines a family of fibers $U_{t_0}\subset U$ of the projection $p:U\to\C^{k+1}$ over the points $C(t_0),\,t_0\in\C$. Choose a point $\alpha\in U\cap p^{-1}(x_0)$ and its metric neighborhood $W\subset \bigcup_{t_0} U_{t_0}$, denote the restriction of $p$ to $W$ by $t$, and its fiber over $C(t_0)$ by $W_{t_0}$.
Lemma \ref{lfub} is applicable to the maps $(W,\alpha)\xrightarrow{(m-m(\alpha),t)}(\C^2,0)\xrightarrow{t}(\C,0)$, as soon as $C\not\subset B_p$ (since $p|_{\Cr p}$ is proper by (II)). This lemma implies:

\vspace{1ex}

the sum of the Milnor numbers of $m|_{W_{t_0}}$ minus the Milnor number of $m|_{W_0}$ at $\alpha$ equals the signed reduced Euler characteristics of $W_{t_0}$ for small $t_0$. 

\vspace{1ex}

IV. Sum up (III) over all $\alpha\in \Cr(m_{U_0})$ (which are finitely many by (II)), and notice that for $t_0\to 0$ all critical points of $m|_{U_{t_0}}$ tend to critical points of $m|_{U_{0}}$, since $p|_{\Cr q}$ is proper by (II). We conclude that the sum of the Milnor numbers of $m|_{U_t}$ minus the same for $m|_{U_0}$ equals:

0 if $x_0\notin B_p$, or

$\mu$ if $x_0$ is a generic point of a weight $\mu$ component of $B_p$.

\vspace{1ex}

By Theorem \ref{thfull}.3 for smaller $n$ and sufficiently generic $m$, the sums of the Milnor numbers in the left hand side of this equality are just the Euler characteristics of $U_{t_0}$ and $U_0$, so the equality proves Theorem \ref{thproj0}.
$\hfill\square$

\begin{rem}
1. We stress that the proof of the new theorem relies on Theorem \ref{thfull}.3 for smaller $n$ and \ref{thfull}.1 for larger $k$, because later the conclusion of the new theorem will be used to prove Theorem \ref{thfull} by induction on increasing $n$ and then decreasing $k$.

2. Incidentally, the characteirstic class of the fiber $\Trop_j U_{t_0}$ equals the intersection prouct of the fans $\ker dp$ and $\Trop_{j+k+1} U$ for all $j$. Once the main results of the paper are proved, this aposteriori follows from the multiplicativity of tropical characteristic classes (Corollary \ref{sledtrop}). However, this observation may be sued for a different proof of the main results.
\end{rem}

{\bf Inductive transversality.} This easy combinatorial fact was used in the preceding proof.
\begin{lemma}\label{lindtransv}
Let $K\subset \R^n$ be a pure $k$-dimensional polyhedral complex submersed by a projection $p:\R^n\to\R^n/L$. Then, for generic hyperplane $H\subset L$, the projection $q:\R^n\to\R^n/H$ submerses $K$ (for $k+\dim L>n$) or even immerses it (for $k+\dim L=n$).
\end{lemma}
\begin{proof}
In case  $k+\dim L>n$, fibers of $p|_K$ have finitely many directions of tangent spaces at their smooth points of the expected dimension $k+\dim L-n>0$. It is enough to choose the hyperplane $H$ transversal to each of these finitely many directions.

In case  $k+\dim L=n$, every maximal dimensional cell $K_i$ of the complex $K$ is mapped injectively by $p$ and thus automatically by $q$. Thus it is enough to show that $k$-dimensional $q(K_i)$ and $q(K_j)$ do not have $k$-dimensional intersection. The proof splits into two cases.

If $K_i$ and $K_j$ have the same affine span, then $p(K_i)\cap p(K_j)=p(K_i\cap K_j)$, and $\dim K_i\cap K_j<k$.

If $K_i$ and $K_j$ have different affine spans, these spans are graphs of functions $\varphi_i$ and $\varphi_j:\R^n/L\to L$, and it is enough to choose $H$ transversal to the (non-trivial) image of $\varphi_j-\varphi_i$.
\end{proof}

{\bf Relation to affine characteristic classes.} The identity $p_*c^{SM}_k U=c^{SM}_k p_*U$ 
proves that the fundamental class of $E_p$ equals $p_*c^{SM}_k U-e_V\cdot c^{SM}_k V$. We need a version for affine classes.
\begin{utver}\label{propdiscrchar}
For 
$U\subset\CC^n\xrightarrow{m}\CC^{k+1}$, we have $m_*\langle U\rangle_{k}=\Trop E_{m|_U}$.
\end{utver}
Recall that the tropical fan $\Trop A$ of a polytope $A\subset\Z^n$ is defined as the corner locus of the support functions of $A$: its $(n-1)$-dimensional cones are dual to the edges of $A$, and their weights equal the lattice lengths of the respective edges.
\begin{proof}
For generic polynomials $g:=(g_1,\ldots,g_k)$ with a given Newton polytope $A\subset\Z^{k+1}$, the Euler characteristics of $U\cap m^{-1}\{g=0\}$ can be computed fiberwise as $e_{gen}\cdot e\{g=0\}-\Trop\{g=0\}\cdot\Trop E_{m|_U}$. We interpret the three Euler characteristics as the 0-dimensional affine characteristic classes of the respective sets, and hence express them in terms of $\langle U\rangle$:

$$\langle\{g=0\}\rangle=\left(\frac{\Trop A}{1-\Trop A}\right)^{k}\;\mbox{ by \cite{kh77}, so }\; e\{g=0\}=-k(\Trop A)^{k+1};$$

$$\langle U\cap m^{-1}\{g=0\}\rangle=\langle U\rangle\cdot m^*\langle\{g=0\}\rangle,\;\mbox{ so }\;$$ $$e(U\cap m^{-1}\{g=0\})=-k \langle U\rangle_{k+1}\cdot m^*(\Trop A)^{k+1} +\langle U\rangle_{k}\cdot m^*(\Trop A)^{k};$$

$$\langle U\cap m^{-1}(y)\rangle=\langle U\rangle\cdot \Trop m^{-1}(y)\mbox{ for generic }y\in V,\; \mbox{ so }\; e_{gen}=\langle U\rangle_{k+1}\cdot m^*\Trop \{y\}.$$
We allow a harmless abuse of notation here, interchangeably treating $(\Trop A)^{k+1}$ as a 0-dimensional tropical fan and its numerical weight, so that in particular $(\Trop A)^{k+1}m^*\Trop \{y\}=m^*(\Trop A)^{k+1}$. Plugging all of this to the right hand side of the initial Euler characteristics identity $\Trop\{g=0\}\cdot\Trop E_{m|_U}=e_{gen}\cdot e\{g=0\}-e(U\cap m^{-1}\{g=0\})$, we get
$$(\Trop A)^k\cdot\Trop E_{m|_U}=\langle U\rangle_{k}\cdot m^*(\Trop A)^{k}=(\Trop A)^{k}\cdot m_*\langle U\rangle_{k}.$$
Since fans of the form $(\Trop A)^k$ generate the space of 1-dimensional fans, and the intersection pairing on the ring of tropical pairs is non-degenerate, we have $\Trop E_{m|_U}=m_*\langle U\rangle_{k}$.
\end{proof}

\section{Relation to Gauss maps}\label{sgauss}

{\bf The nondegeneracy.} We show that the rank drop locus $Z(W)$ is reduced for a generic form $w$.
\begin{utver}\label{propnondeg1}
1. For a $q$-dimensional variety $U\subset\CC^n$ and a generic invariant form $w=w_0\wedge\ldots\wedge w_k$, the set $Z(w)\subset\sm U$ is $k$-dimensional or empty.

2. If it is not empty, then, for generic $w=w_0\wedge\ldots\wedge w_k$ and a generic point $x\in Z(w)$, represent $w_i$ as $dm_i$ near $x$. If $x$ is an isolated solution of the equations $$m_0=\cdots=m_k=0$$ on $Z(w)$, then, on $\sm U$, these equations define a reduced isolated complete intersection singularity (ICIS), and this $(x,w)$-ICIS is nondegenerate (i.e. has Milnor number 1).
\end{utver}
For the proof, consider the {\it tautological space} with its projections $u$ and $\pi$ to the multipliers: $$S:=\{(x,L)\,|\, L+T_xU\subsetneq T_x\CC^n\}\subset\sm U\times Gr(n,n-k-1).$$ We identify every tangent space $T_x\CC^n$ with $\C^n:=T_{(1,\ldots,1)}\CC^n$ here, since $\CC^n$ is a group.

Since $u$ is a locally trivial fibration into Schubert varieties, $S$ is irreducible. Fibers of $\pi$ are exactly the sets $Z(w)$ once we identify the form $w$ with the space $\ker w\subset T_{(1,\ldots,1)}\CC^n$. 

{\it Proof of Proposition \ref{propnondeg1}.1.} From the fibration structure of $u$, we see $\dim S=n-1+k(n-k-1)$. Thus the generic fiber of $\pi$ has dimension $k$ (if $\pi$ is dominant).$\hfill\square$

For part 2, recall some properties of the projectively dual variety $H^\vee$ to an analytic hypersurface $H:=\{h=0\}\subset\C^m$ near its smooth point $z$. If $\dim H^\vee=r$, we have the following.

i) For every $y\in H$, the tangent hyperplane $T_{y} H$ is tangent to $H$ at every point of some $(m-r-1)$-dimensional plane $R_{y}$.

ii) For generic $y\in H$, the quadratic form $d^2h|_{T_{y} H}$ has $(m-r-1)$-dimensional kernel (thus equal to $R_{y}$).

The proof follows from comparing $d^2h|_{T_{y} H}$ to the differential of the Gauss map of $H$.

From this, we can deduce the following logarithmic version of Proposition \ref{propnondeg1}.2.
\begin{lemma}\label{lgauss}
Given an analytic manifold $V\subset\C^n$ and a linear map $\mu:\C^n\mapsto\C^k$, assume 

-- the critical locus $Z(\mu)$ of $\mu|_V$ has the expected dimension $(n-k-1)$ at its point $z$, 

-- the map $\mu$ is finite on it near $z$, and

-- $\dim (T_z V+\ker\mu)=n-1$. 

Then there exist arbitrarily close $\tilde\mu$ and $\tilde z$, such that the equations $\tilde\mu=0$ define on $V$ a nondegenerate ICIS at $\tilde z$.
\end{lemma}
\begin{proof}
First, assume $\dim V=n-1$. Applying $(i)$ with $H=V$ we see that $r\geq n-k-1$. Thus, by (ii), the restriction of $h$ onto $\tilde\mu=0$ is a nondegenerate hypersurface singularity. Then the restriction of $\tilde\mu$ onto $h=0$ is a nondegenerate ICIS (because these two singularities have homeomorphic Milnor numbers: if one is a homotopy sphere, so is the other one).

If $q:=\dim V<n-1$, the equality $\dim (T_z V+\ker\mu)=n-1$ allows us to factor $\mu$ as $\C^n\xrightarrow{\nu }C^{q+1}\xrightarrow{\mu'}\C^k$ so that $T_z V\cap \ker\nu=\{0\}$. Then the statement of lemma for $V$ reduces to the one for the smooth hypersurface $H:=\nu(V)$, which we have already proved.
\end{proof}

{\it Proof of Proposition \ref{propnondeg1}.2.} The dimension count shows that, for generic $(x,w)\in S$, we have $\dim(T_x U+\ker w)=n-1$. This allows to apply Lemma \ref{lgauss} to the preimage $(V,z)$ of $(U,x)$ under the map $\exp:\C^n\to\CC^n$ (choosing a branch of the logarithm near $x$). $\hfill\square$

\vspace{1ex}

{\bf Restrictions of equidimensional maps.} In order to further study the projection $\pi:S\to Gr(n,n-k-1)$, we recall some generalities about arbitrary equidimensional maps.

Let $\pi$ map an analytic space $A$ of (pure) dimension $a$ to a smooth variety $B$ of dimension $b$, so that every fiber has the expected dimension $a-b$. We allow $\pi$ to be non-proper and some of its fibers to have components in $\sing A$. Recall that the Bertini discriminant set $B_\pi\subset B$ is the closure of the union of $\pi\sing A$ and critical values of $\pi|_{\sm A}$.

\begin{utver}\label{pfib0}
For every $x\in A$, its fiber has a structure of a cycle (denoted by $\pi^* x$), which is natural in the following sense: 

1. The cycle $\pi^* x$ is reduced for generic $x$, particularly for $x\notin B_\pi$.

2. For any $x$, we have $\pi^* x=\lim_{\tilde x\to x}\pi^*\tilde x$ (in the sense of Definition \ref{defconserv}).
\end{utver}
Recall how to explicitly define the weight of a component $C\subset\pi^{-1}x$: choose a germ $H$ of a generic $b$-dimensional plane section of $A$ at a generic point $y\in C$, and define the weight of $C$ in $\pi^* x$ as the degree of the finite map $\pi|_H:H\to A$ at $x$. 

\begin{utver}\label{pfib1}
Let $B'\subset B$ be smooth, and denote  the restriction of $\pi$ to its preimage $A':=\pi^{-1}B'$ by $\pi'$.

1. $B_{\pi'}\subset B_\pi$.

2. $A'$ has pure dimension.

3. If $B'\not\subset B_\pi$, then $\pi'^* x=\pi^* x$ for $x\in B'$.
\end{utver}
Part (1) is tautological. Part (2) is not, but follows by induction on codimension of $B'\subset B$. Then part (3) follows for generic $x\in B'$ because both cycles are reduced (Proposition \ref{pfib0}.1), and then for any $x$ by conservation of intersection numbers (Proposition \ref{pfib0}.2).

\vspace{1ex}

We apply this to the map $\pi:S\to G,\,G:=Gr(n,n-k-1)$ from the preceding section. As before, every rational plane $P\in G$ is identified with 

-- an invariant $(k+1)$-form $w$ on the torus $\CC^n$, such that $\ker w=P$ in $T_1\CC^n$, and 

-- the projection $p:\CC^n\to\CC^{k+1}$ such that $\ker dp=P$. 

Furthermore, the form $w$ gives rise to the set $Z(w)\subset\sm U$, and a cycle $\Cr p$ supported on it (though $p$ is defined only for rational planes $P$, the respective cycle can be actually defined in general, as in Remark \ref{rembert}.1.) However, the projection $\pi$ defines one more cycle $\pi^* P$, whose image under the projection $S\to\sm U$ has the same support $Z(w)$.

We want to show these two cycles are equal: at a generic point of $Z(w)$, the Milnor number of $p|_U$ equals the number of branches of $Z(\tilde w)$ into which the set $Z(w)$ slits as we perturb $w$ to a close generic form $\tilde w$.
Expectedly, if $U$ is a hypersurface, the equality boils down to Milnor's theorem that the Milnor number equals the degree of the gradient map. 

For the general case and subsequent applications, it will be convenient to consider a pencil of planes $P_t\in G$, denoting respective forms and projection maps by $w_t$ and $p_t$.
\begin{utver} \label{ppencil}
1. Given a rational $P\in G$ such that the set $Z(w)$ is $k$-dimensional, and the projection $p$ is finite on it, we have the equality of cycles $\pi^*P=\Cr p$ supported at $Z(w)$.

2. For any pencil $P_t$ parameterized near $P_0=P$ by $t$ in a small disc $B\subset\C$, the set $\bigcup_{t\in B} Z(w_t)$ has pure dimension, and we have $\Cr p_0=\lim_{t\to 0} \Cr p_t$.

3. For a generic such pencil, the cycles $\Cr p_t=\pi^*P_t$ are reduced for non-zero $t\in B$.
\end{utver}
\begin{proof}
Part 2: pure dimension follows by Proposition \ref{pfib1}, and the equality by Corollary \ref{conserv}.

Part 3 follows by Proposition \ref{propnondeg1}.

Part 1:  
choosing a generic pencil $P_t$ with $P_0=P$, we see $$\Cr p_0=\lim_{t\to 0} \Cr p_t=\lim_{t\to 0} \pi^*P_t=\pi^*P_0.$$
Here the first equality is Part 2, the second is because both $\Cr p_t$ and $\pi^*P_t$ are reduced for $t\ne 0$ (by Propositions \ref{pfib0} and \ref{propnondeg1} respectively), and the third equality by Proposition \ref{pfib0}.
\end{proof}

\vspace{1ex}

{\bf Proof of Theorem \ref{thclassif}.} This is the last ingredient to prove the main results.
\begin{lemma}
Every point $x\in V$ belongs to a nontrivial torus coset contained in $V$.
\end{lemma}
\begin{proof}
Choose $\tilde x\in\C^n$ so that $\exp(\tilde x)=x$, and denote the germ of the preimage $\exp^{-1}(X)$ at $\tilde x$ by $\tilde V$. The assumption of the theorem implies that $\tilde V$ is dual defective (i.e. its projectively dual $\tilde V^\vee$ has codimension higher than 1), thus it is ruled by affine planes of positive dimension (which are fibers of the projectivized normal bundle of $\tilde V^\vee$, since $\tilde V^{\vee\vee}=\tilde V$). In particular, $\tilde V$ contains a line $\tilde x+\R\cdot v$, thus $V$ contains a non-trivial torus coset $x\cdot\overline{\exp(\C\cdot v)}$.
\end{proof}
For every torus $T\subset\CC^n$, define the Zariski closed set $V_T\subset V$ as the union of all cosets of $T$ contained in $V$. Lemma implies that $V=\bigcup_T V_T$ over all non-trivial tori $T$. Since they are contably many, we have $V=V_T$ for some $T$. This is the statement of Theorem \ref{thclassif}.1.

Theorem \ref{thclassif}.2 then follows by induction on $r$, because $\Trop_{k+1} \C^*\times V=\Trop \C^*\times\Trop_k V$.

\vspace{1ex}

{\bf Proof of Theorem \ref{thfull}.} 
I. If $k=0$, then Theorem \ref{thproj0} gives the equality of 0-dimensional cycles $B_p=E_p$, and passing to their degrees we have $\deg\Cr p=\deg B_p=\deg E_p=e(U)$ (the latter equality is from counting $e(U)$ fiberwise). This proves part (3), and implies that $\Trop\Cr p$ is the point $0$ with the same weight $\deg\Cr p=e(U)$ for all projections $p$ satisfying the genericity assumption of the theorem. This tropical fan is by definition $\Trop_0 U$, and its weight is $e(U)=\langle U\rangle_n$: this proves parts (1) and (2).

\vspace{1ex}

II. If $k>0$, but the projection $\pi$ (defined after Proposition \ref{propnondeg1}) is not dominant, then, by  Theorem \ref{thclassif}, we have $U=U'\times\C^*$. Apply Theorem \ref{thfull} to $U'\subset\CC^{n-1}$ by induction on $n$.

\vspace{1ex}

We now consider the remaining case: $\dim\sing U<k\ne 0$ and projection $\pi$ is dominant. 

\vspace{1ex}

III. Assume a $(k+1)$-form $w$ is rational (so that the cycle $Z(w)$ equals $\Cr p$ for the respective projection $p$) and generic (in the sense that the fan $\Trop Z(w)$ equals its generic instance $\Trop_k U$). Then we see from Theorem \ref{thproj0} that $$p_* Z(w)=p_* \Cr p=B_p=E_p.$$
Passing to tropicalizations, we see $p_*\Trop_k U=\Trop E_p=p_* \langle U\rangle_k$ by Proposition \ref{propdiscrchar}. Since this holds for a generic projection $p$, we conclude $\Trop_k U= \langle U\rangle_k$. Parts 1 and 2 are proved.

\vspace{1ex}

IV. It remains to prove that $\Trop_k U=\Trop Z(w)$ once the plane $\ker w$ immerses $\Trop_k U$ and submerses $\Trop_j U$ for $j>k$. Choose a generic rational 1-dimensional pencil of forms $w_t$ with $w_0=w$, in the sense that $\Trop Z(w_t)=\Trop_k U$ for small $t\ne 0$. For such $t$, we see $$p_*\Trop Z(w_0)=\Trop B_p=\Trop E_p=p_*\langle U\rangle_k=p_*\Trop_k U=p_*\Trop Z(w_t).$$ Here the equality $B_p=E_p$ is again by Theorem \ref{thproj0}, and the second last one by (III).

Thus we can apply the tropical isotopy theorem \ref{tropisot}, having checked that 

-- the set $\bigcup_t Z(w_t)$ has pure dimension by Proposition \ref{ppencil};

-- $Z(w_0)=\lim_{t\to 0} Z(w_t)$ by the same proposition;

-- the plane $\ker dp$ submerses or immerses the relevant fans $\Trop Z(w_t)=\Trop_k U$. 

This gives $\Trop Z(w_0)=\Trop Z(w_t)=\Trop_k U$.

\section{The proof of the tropical isotopy theorem}\label{sisot}

{\bf Nice toric morphisms.} Given a $k$-dimensional subvariety $V\subset\CC^n\xrightarrow{p}\C^r$, we want to construct a nice toric compactification $X_\Sigma\supset\CC^n$ of $V$, on which we could nicely extend the coordinate projection $p$. More specifically:

1) A nice toric compactification means {\it tropical} (\cite{tev}). Geometrically, this means that the closure $\bar V\subset\CC^n$ does not meet toric orbits of codimension higher than $k$. Combinatorially, the toric fan $\Sigma$ refines the tropical fan $\Trop V$ (i.e. $\Trop V$ is in the $k$-skeleton of $\Sigma$).

2) A nice toric morphism means {\it weakly semistable} (\cite{abrkaru00}). Geometrically, this means that the map of smooth complete toric varieites $p:X_\Sigma\to X_{\Sigma'}$, extending the projection $X_\Sigma\supset\CC^n\xrightarrow{p}\CC^r\subset X_{\Sigma'}$, is flat with reduced fibers. Combinatorially, the image of every cone $C\in\Sigma$ is some cone $C'\in\Sigma'$, and, moreover, $p_*(C\cap\Z^n)=C'\cap\Z^r$.

We describe a weakly semistable morphism in the smallest interesting case $p:\CC^3\to\CC^2,\, p(x,y,z)\mapsto(x,y)$, the only one that we need for the tropical isotopy lemma. 
\begin{utver}\label{propws} 1. A pair of fans $\Sigma$ in $\R^3$ and $\Sigma'$ in $\R^2$ describing a weakly semistable morphism has the following structure:

- the consecutive rays of $\Sigma'$ are generated by primitive vectors $v_1,\ldots,v_N=v_0\in\Z^2$, such that $\det(v_{i-1},v_i)=1$ for $i=1,\ldots,N$;

-- the rays of $\Sigma$ are generated by vectors $u_{\pm}:=(0,0,\pm 1)$ and $u_{i,j}:=(v_i,j),\,j=a_i,\ldots,b_i$;

-- the combinatorial structure of $\Sigma$ is as on the picture below.
\begin{figure}
\begin{center}
\includegraphics[scale=0.35]{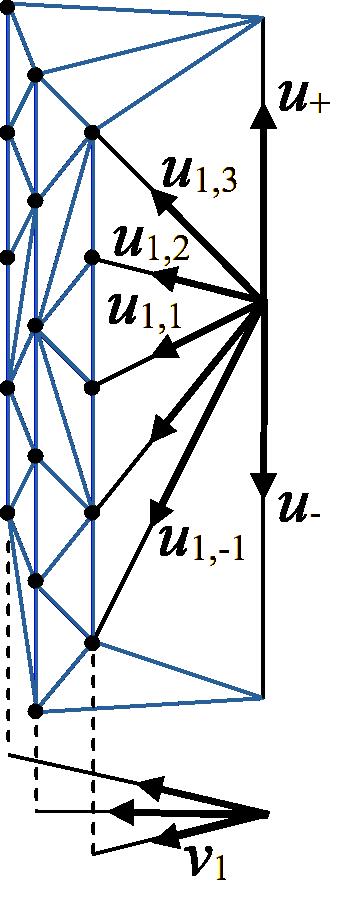}
\end{center}
\end{figure}

2. For any rational fan $S$ in $R^3$, there exists a sufficiently divisible $m\in\Z$ such that the image of $S$ under the dilatation $d_m(x,y,z):=(x,y,mz)$ refines to a complete fan $\Sigma$, which defines a weakly semistable toric morphism extending $p$.
\end{utver}
\begin{proof} Part 1 is clear from the combinatorial definition of weak semistability. For part 2:

-- take $\Sigma'$ to be a smooth fan refining $p(1$-skeleton of $S)$;

-- let $S_1$ be the 1-skeleton of $(2$-skeleton of $S)\cap p^{-1}(1$-skeleton of $\Sigma')$;

-- choose $m$ to be GCD of $x_i$ and $y_i$ over all lattice vectors $(x_i,y,_i,z_i)$ generating the rays of $S_1$, then the sought $\Sigma$ can be chosen to have 1-skeleton $d_m(S_1)$.
\end{proof}
\begin{sledst}
Given a surface $\{f(x,y,z)=0\}\subset\CC^3\xrightarrow{p}\C^2$, one can extend the projection $p(x,y,z)=(x,y)$ to a weakly semistable morphism of toric varieties $p:X_\Sigma\to X_{\Sigma'}$, such that $X_\Sigma$ is the tropical compactification for the surface $\{f(x^m,y^m,z)=0\}$ for some (sufficiently large and divisible) $m\in\Z$.
\end{sledst}
\begin{proof}
Apply the preceding proposition to $S=\Trop\{f=0\}$.
\end{proof}
\begin{rem}
For our purpose, it would be actually enough to extend the projection $p$ to any (not necessarily weakly semistable) morphism of smooth complete toric varieties, such that the domain is a tropical compactification of a given surface. However, the simplest known way to construct such a morphism is as sketched in the proof of Proposition \ref{propws}, and the result happens to be weakly semistable for the same price.
\end{rem}

{\bf Tropicalizations with components at infinity.}  
We recast (in an obvious way) some elements of tropical geometry for very affine varieties having ``components at infinity''.
\begin{defin}\label{deftropinf}
Let $Z\subset X_\Sigma$ be  an algebraic set (or cycle) of pure dimension $k$ in a smooth complete toric variety, 
then its {\it tropical fan} is the $k$-dimensional skeleton of the toric fan $\Sigma$, with the weight of each $k$-dimensional cone $C\in\Sigma$ defined as the intersection index $Z\cdot\overline{(C\mbox{-orbit})}$; in terms of \cite{fs}, it is the Minkowski weight of the fundamental cycle of $Z$.
\end{defin}
Note that the intersection index with a compact subvariety (such as the closure of an orbit) is  defined (e.g. via cohomology multiplication) even if the intersection has positive dimension.
\begin{utver}\label{tropprops}
1 -- \cite{fs}. The intersection index of sets $Z$ and $Z'$ of complementary dimensions in $X_\Sigma$ equals the intersection number $\Trop Z\cdot\Trop Z'$.

2 -- \cite{st}. Given algebraic cycles $Z$ and $Z'$ in the domain and range of a weakly semistable morphism $p:X_\Sigma\to X_{\Sigma'}$, we have $\Trop p^{-1} Z'=p^* \Trop Z'$ and $\Trop p_* Z=p_* \Trop Z$.
\end{utver}
Here $p_* V$ for a variety $V$ is $($degree of $p|_V)\cdot V$ if $\dim p(V)=\dim V$ and 0 otherwise. 

\begin{utver}\label{utvpos}
1. If $H$ is a positively weighted $k$-dimensional tropical fan in the $k$-skeleton of a toric fan $\Sigma$, and $Z$ is a positively weighted codimension $k$ algebraic cycle in the smooth projective $n$-dimensional toric variety $X_\Sigma$, then $H\cdot\Trop Z\geq 0$.

2. There exists a $k$-dimensional subvariety $V\subset\CC^n\subset X_\Sigma$ such that $\Trop V\subset(k$-skeleton of $\Sigma)$ and $\Trop V\cdot\Trop Z>0$.
\end{utver}
This is not a tautology, because $\Trop Z$ may have negative weights, and $H\cdot\Trop Z$ may be negative unless the $k$-skeleton of $\Sigma$ contains $H$. We need only the following special case.

{\it Proof under the additional assumption that $H$ is the tropical fan of an algebraic variety $W\subset\CC^n$.} For Part 1, $\Trop W\subset(k$-skeleton of $\Sigma)$ implies that $W$ is tropically compactified by $X_\Sigma$. Thus $Z\cap\bar gW\subset X_\Sigma$ is 0-dimensional for generic $g\in\CC^n$, so $\Trop Z\cdot H=Z\cdot\bar gW >0$. For Part 2, take $\bar V$ to be a generic $k$-codimensional plane section of $X_\Sigma$. $\hfill\square$

If $Z$ is the closure of a subvariety of a torus $\CC^n$ in its tropical compactification $X_\Sigma$, Definition \ref{deftropinf} gives the usual notion of a tropical fan. We however arrive at our higher generality when studying tropical fans of a family of varieties $V_t\subset\CC^n$ as in Theorem \ref{tropisot}: the members $V_t$ for all small $t\ne 0$ have the same tropical fan denoted by $\Trop V_\bullet$, which in general differs from $\Trop V_0$. 
To analyze the difference, it is convenient to choose a toric fan $\Sigma$ refining both $\Trop V_\bullet$ and $\Trop V_0$, so that the respective toric variety $X_\Sigma$ is smooth and projective. In it, we can take the closures $\bar V_t\subset X_\Sigma$ and consider their limit $\lim_{t\to 0}\bar V_t\subset  X_\Sigma$ as an algebraic cycle defined as the divisor $t=0$ on $\bar V\subset X_\Sigma\times\C$ (where the function $t$ is the projection to $\C$, and the components of the divisor come with natural weights measuring the order of zero locus).

Note that, besides $\bar V_0$, this limit $\lim_{t\to 0}\bar V_t$ may have {\it boundary components} $Z_i\subset X_\Sigma\setminus\CC^n$, and the following observation is obvious from the conservation of intersection indices.
\begin{utver}\label{propconserv}
We have $\Trop V_\bullet=\Trop \bar V_\bullet=\Trop\lim_{t\to 0}\bar V_t=\Trop V_0+\sum_i\Trop Z_i.$
\end{utver}
It is useful to look at the fans $\Trop Z_i$ when analysing the difference $\Trop V_\bullet-\Trop V_0$, though usually they have cones with negative weights. On a positive side, we have the following.
\begin{utver}\label{semicont}
1. For any positively weighted tropical fan $H$ of complimentary dimension, we have $H\cdot\Trop V_0\leqslant H\cdot\Trop V_\bullet.$

2. If $\lim_{t\to 0}\bar V_t$ has a boundary component $Z$, then the inequality is strict for some $H$.
\end{utver}
Note that, in contrast to Proposition \ref{utvpos}.1, we do not assume $H\subset(k$-skeleton of $\Sigma)$ here.

{\it Proof under the additional assumption that $H$ is the tropical fan of an algebraic variety $W\subset\CC^n$.} Pick $g\in\CC^n$ such that $(gH)\cdot V_t=H\cdot\Trop V_t$ for all small $t$ (including 0), and notice that $(gH)\cdot V_0\leqslant (gH)\cdot V_t$ for small $t\ne 0$ by semicontinuity of the intersection number.
For Part 2, multiply the equality of Proposition \ref{propconserv} by the fan $\Trop V$ of Proposition \ref{utvpos}.2.
$\hfill\square$

\vspace{1ex}

{\bf Tropical compactifications of families.} The key to the general case (which we do not elaborate because do not use it here) is the following fact of independent importance.
\begin{utver}\label{familycomp}
For every family of varieties $V_t\subset\CC^n$ as in Theorem \ref{tropisot}, there exists a fan $\Sigma$ such that the projective smooth toric variety $X_\Sigma$ is a {\it codimension $1$ tropical compactification} of $\lim_{t\to 0}\bar V_t\subset X_\Sigma$ in the following sense. 
\end{utver}
\begin{defin}
1. The {\it ambient orbit} of a subvariety $V\subset X_\Sigma$ is the the maximal orbit of the toric variety $X_\Sigma$ intersecting $V$ (note that its closure contains $V$).

2. $V$ is said to be {\it in codimension $r$} of $X_\Sigma$, if its ambient orbit has codimension at least $r$.

3. $V$ is said to be tropically compactified by $X_\Sigma$, if it is tropically compactified by the closure $Y$ of its ambient orbit (i.e. if it evades orbits of $Y$ of codimension higher than $\dim V$).

4. A toric variety $X_\Sigma$ is a codimension $r$ tropical compactification of its algebraic cycle $Z$, if every component of $Z$ is in codimension $r$ of $X_\Sigma$ and is tropically compactified by it.
\end{defin}
\begin{exa}\label{exatc}
A toric threefold $X$ tropically compactifies a curve $C\subset X$ if no component of $C$  touches 0-orbits, and the components touching the 3-orbit do not touch 1-orbits. In particular $X$ is then automatically a codimension 1 tropical compactification of $C$ (i.e. no component of $C$ is the closure of a 1-orbit).
\end{exa}
{\it Proof of Proposition \ref{familycomp}.}
Let $S$ be the projection of $\Trop V\subset\R^n\oplus\R^1$ to $\R^n$, and let $\Sigma$ be a smooth complete fan refining $S$. Since $S$ is, equivalently, the tropical fan of the projection $W$ of the unfolding $V\subset\CC^n\times\CC^1$ to $\CC^n$ (\cite{st}), the toric variety $X_\Sigma$ is a tropical compactification for $W$. Since $\bigcup_i Z_i\subset\overline{W}\setminus\CC^n\subset X_\Sigma$, we conclude that $Z_i$ evade orbits of codimension higher than $\dim W=\dim V=k+1$.
On the other hand, $\dim Z_i=k$, so:

1) if the maximal orbit $O$ intersecting $Z_i$ has dimension $n-1$, then $\bar O\supset Z_i$ is a tropical compactification, as $Z_i$ evades ordbits of codimension higher than $k$ in the toric variety $\bar O$ (because they have codimension higher than $k+1$ in $X_\Sigma$);

2) the maximal orbit $O$ intersecting $Z_i$ cannot have dimension below $n-1$, because otherwise $Z_i$ would have to intersect at least one codimension $k$ orbit in the toric variety $\bar O\supset Z_i$, which would have codimension higher than $k+1$ in $X_\Sigma$.
$\hfill\square$
\begin{rem}
1. The interaction of pull-backs and push-forwards with taking $\Trop$ is more subtle than Proposition \ref{tropprops}.2 for non weakly semistable morphisms. However this is easy to repair for cycles that are tropically compactified in a small codimension. 

2. For toric varieties $Y\subset X$ and a subvariety $V\subset Y$, the tropical fan $\Trop V$ with respect to $X$ is easily reconstructed from the one with respect to $Y$. 

We do not elaborate these two points here, because do not use them.
\end{rem}

{\bf Proof of Theorem \ref{tropisot}} proceeds in three steps:  prove the simplest case $(n,r,k)=(3,2,1)$, then observe that extending this to arbitrary $(n,r)$ and then $k$ is mostly a general nonsense. 
\begin{lemma}
The theorem holds true for $(n,r,k)=(3,2,1)$, i.e. for a family of curves $V_t$ in the domain of the projection $p:\CC^3\to\CC^2,\,p(x,y,z)=(x,y)$. 
\end{lemma}
To sketch a proof (neglecting the rigour), look at the boundary components of $\lim_{t\to 0}\bar V_t$ in a nice toric compactification, and notice three possibilities:

-- if there is a ``vertical'' boundary component $Z$ (mapped to a point by $p$), then $p_*$ is not injective on $\Trop V_\bullet$: it contains two rays of $\Trop Z$ with the same projection;

-- if there is a ``horizontal'' boundary component (mapped to a curve $W$ by $p$), then the difference $p_*\Trop V_\bullet-p_*\Trop V_0$ is nontrivial: it is estimated from below by $\Trop W$;

-- finally, the absence of any boundary components proves the lemma by Proposition \ref{propconserv}.
\begin{proof}
Choose a toric compactification $X_{\tilde\Sigma}\supset\CC^3$ as in Proposition \ref{familycomp}, and apply Proposition \ref{propws}.2 to the union of the tropical fans of the components of the curve $\lim_{t\to 0}\bar V_t\subset X_{\tilde\Sigma}$. We obtain 

\begin{center}
\begin{tikzcd}
V_t \arrow[d, phantom, sloped, "\subset"]   & U_t \arrow[l, maps to] \arrow[d, phantom, sloped, "\subset"]                & \lim\limits_{t\to 0}\bar U_t \arrow[d, phantom, sloped, "\subset"] \\
\CC^3 \arrow[d] & \CC^3 \arrow[d] \arrow[l, "\varphi"', maps to] \arrow[r, phantom, sloped, "\subset"] & X_{\Sigma} \arrow[d, "p"]              \\
\CC^2           & \CC^2 \arrow[l, "\varphi'"', maps to] \arrow[r, phantom, sloped, "\subset"]          & X_{\Sigma'} 
\end{tikzcd}
\end{center}
where $\varphi'$ is a sufficiently divisible covering of tori, $U_t$ is the preimage of $V_t$ under the induced covering $\varphi$, and $p$ is weakly semistable as in Proposition \ref{propws}.1. The equalities $\Trop \bar U_t=\varphi^* \Trop \bar V_t$ and $\Trop \lim_{t\to 0}\bar U_t=\varphi^* \Trop \lim_{t\to 0}\bar V_t$  
have two important implications:

1) it is enough to prove the lemma for the family $U_t$;

2) the choice of $X_{\tilde\Sigma}$ assures that the toric variety $X_{\Sigma}$  tropically compactifies the curve $\lim_{t\to 0}\bar U_t$ as in Example \ref{exatc}. Look at its boundary component $Z$.

If $p_* Z\ne 0$, then $p_* Z$ is a boundary component of the family $p_* V_t$, 
thus $$p_*\Trop U_\bullet=\Trop p_* U_\bullet>\Trop p_* U_0=p_*\Trop U_0$$
(with $=$ and $>$ by Propositions \ref{tropprops}.2 and \ref{semicont}.2), which contradicts $p_*\Trop U_\bullet=p_*\Trop U_0$.

So $p_* Z=0$ for every boundary component $Z$. From Example \ref{exatc} in the notation of Proposition \ref{propws}.1, the curve $Z$ belongs to the $u_{i,j}$-orbit of the toric variety $X_\Sigma$ for some $i$ and $j\in[a_i,b_i]$, and $\Trop Z$ consists of three rays: $\R_+\cdot u_{i,j}$ with a negative weight and $\R_+\cdot u_{i,j\pm 1}$ with positive weights (here we set $u_{i,j\pm 1}:=u_\pm$ if $j\pm 1\notin[a_i,b_i]$).

An obivious but crucial combinatorial observation is that any nontrivial linear combination of such tripods has a positively weighted ray $\R_+\cdot u_\pm$ or positively weighted rays $\R_+\cdot u_{i,j}$ and $\R_+\cdot u_{i,j'}$ for $j\ne j'$. Thus such a ray or pair sits in $\Trop U_\bullet$ by Proposition \ref{propconserv}, if only $\lim_{t\to 0}\bar U_t$ has a boundary component.
Since such a ray or pair would destroy injectivity of $p_*$ on $\Trop U_\bullet$, there are no boundary components, and $\Trop U_\bullet=\Trop U_0$ by the same Proposition \ref{propconserv}.
\end{proof}
We continue the proof of Theorem \ref{tropisot}, proceeding from the trick to the general nonsence.
\begin{lemma}
For $k=1$, the theorem reduces to the known case $(n,r)=(3,2)$.
\end{lemma}
\begin{proof}
For a projection $q:\CC^n\to\CC^3$, defined as generic $(\varphi,\psi)$ in the diagram
\begin{center}
\begin{tikzcd}
\CC^2\times\CC^1 \arrow[d, "p'"] & \CC^r\times\CC^{n-r} 
\arrow[d, "p"] \arrow[l, "(\varphi{,}\,\psi)"', maps to] \\
\CC^2           & \CC^r \arrow[l, "\varphi"', maps to] 
\end{tikzcd}
\end{center}
the assumptions of the theorem are satisfied for the family $q_*V_t$ and its projection $p'$. Thus $q_*\Trop V_\bullet=\Trop q_* V_\bullet=\Trop q_* V_0=q_*\Trop V_0$. Since this fan has no rays in $0\oplus\R^2\subset\R^1\oplus\R^2=\R^3$ (by injectivity of $p'_*$ on it), the fans $\Trop V_\bullet$ and $\Trop V_0$ have no rays in 
$0\oplus\R^{n-r}\subset\R^r\oplus\R^{n-r}=\R^n$. Since moreover any two rays outside this subspace are distinguished by a generic projection $q=(\varphi,\psi)$, the equality $q_*\Trop V_\bullet=q_*\Trop V_0$ for almost all $q$ implies the sought $\Trop V_\bullet=\Trop V_0$.
\end{proof}
\begin{lemma}
The case of arbitrary $k$ in Theorem \ref{tropisot} reduces to $k=1$.
\end{lemma}
\begin{proof}
Assume towards the contradiction that the theorem is proved for $k=1$ but fails for some family of varieties $V_t$ of higher dimension $k$. Choose generic lattice polytopes $B_1,\ldots,B_{k-1}$ in the sense that the intersection product $F:=\prod_i\Trop B_i$ intersects $p_* S$ only at $0$. Then $p_*$ is injective on the 1-dimensional fan $\Trop V_\bullet\cdot p^*F$: injectivity over $0\in\R^r$ follows $\Trop V_\bullet\cap p^{-1}(0)=\{0\}$ (otherwise $\Trop V_\bullet$ had a cone on which $p_*$ is not injective, then it had a $k$-dimensional cone with this property, and this would contradict the injectivity of $p_*$ on an open dense subset of $\Trop V_\bullet$). Injectiviety over other points of $\R^r$ follows because these points do not belong to $p_* S$ by the choice of $B_i$'s.

Define the Laurent polynomial $g_i$ as a generic linear combination of monomials from $B_i$, then both the generic and the central element in the family of curves $W_t:=V_t\cap\{g_1=\cdots=g_{k-1}=0\}$ has the tropical fans $\Trop W_\bullet=\Trop V_\bullet\cdot p^*F$ and $\Trop W_0=\Trop V_0\cdot p^*F$ respectively, by the tropical correspondence for the equivariant Chow ring of the torus $\CC^n$. Furthermore:
$$p_*\Trop W_\bullet=p_*(\Trop V_\bullet\cdot p^*F)=F\cdot p_*\Trop V_\bullet=F\cdot p_*\Trop V_0=p_*(\Trop V_0\cdot p^*F)=p_*\Trop W_0.$$

To summarize, we have proved that if a family of varieties $V_t$ satisfies the assumption of the tropical isotopy lemma, then the family of curves $W_t$ does. Since we have tropical isotopy in dimension 1, we conclude for ``almost all'' fans $F$: 
$$\Trop V_\bullet\cdot p^*F=\Trop W_\bullet=\Trop W_0=\Trop V_0\cdot p^*F.$$

This implies that the fan $G:=\Trop V_0-\Trop V_\bullet$ as a set projects to a set of dimension smaller than $k$, i.e. none of its $k$-dimensional cones projects injectively. Since $\Trop V_\bullet$ cannot have such cones, all cones in the fan $G$ have non-negative multiplicities.

Assuming (towards contradiction) that $G\ne 0$, notice that 
$G\cdot H>0$ for $H:=\Trop($generic codimension $k$ plane section of $X_\Sigma)$. 
Rewriting this by definition of $G$ as $(\Trop V_0-\Trop V_\bullet)\cdot H>0$, we get the contradiction with Proposition \ref{semicont}. 
Note that $H$ that we use 
satisfies the additional assumption under which we have proved Proposition \ref{semicont}.
\end{proof}

\section{Discussion}\label{sdiscus}

{\bf Other reductive groups.} The {\it equivariant product} $U\cdot V$ of varieties $U$ and $V\subset\CC^n$ is the intersection of $U$ and the shifted copy $g\cdot V$ 
for generic $g\in\CC^n$. It depends on $g$ as an algebraic set, but:
i) for $\dim U+\dim V=n$, its cardinality $|U\cdot V|$ does not depend on $g$;

ii) for any dimension, it does not depend on $g$ as an element of the {\it equivariant Chow group}
$$C^k:=\{\sum\nolimits_i a_iU_i\,:\, a_i\in\Z,\, \codim U_i=k\}/\{\sum\nolimits_i a_iU_i\,:\, \sum\nolimits_i a_i|U_i\cdot V|=0 \mbox{ for all }k\mbox{-dimensional } V\}.$$
It thus induces a ring structure on $C:=C^0\oplus\cdots\oplus C^n$, called the {\it equivariant Chow ring} of $\CC^n$.

The map $U\to\Trop U$ descends to the isomoprhism $C\to($ring of tropical fans$)$, so we interchangeably refer to elements of $C^k$ and respective codimension $k$ tropical fans.

\vspace{1ex}

A version of Theorem \ref{th1} for varieties $U\subset\C^n$ (with less constructive genericity assumptions) dates back at least to \cite{piene} and \cite{schuerm},
but it is hard to imagine how the rest of our work could be transferred to $\C^n$. 
A fundamental reason is that the group $\C^n$ has no equivariant Chow ring (because (ii) fails). That is why the computations in e.g. \cite{schuerm} are done in the homology ring of the variety $U$ itself, which is by far not as rich as the equivariant Chow ring, and thus cannot keep enough information for inductive constructions like \ref{constr1}.

\vspace{1ex}

On the other hand, one can think of extending our results to subvarieties in arbitrary reductive algebraic groups $G$, see e.g. \cite{valya} for the case $U=G$. How to extend our results to an arbitrary smooth subvariety $U$ in a spherical variety?

\vspace{1ex}

{\bf Arbitrarily singular varieties.} 1. How to extend our results to the cycle $Z(w_0\wedge\ldots\wedge w_k)$ the fan $\Trop_k U$ and the CSM class $c^{SM}_k(U)$ for $k=\dim\sing U$? 

\vspace{1ex}

While the weights of the cycle can still be defined as in Remark \ref{rembert}.1, and Theorem \ref{thproj0} admits the case $k=\dim\sing U$, some other steps in our work become subtle in this case. Note that, for a sufficiently generic form $w:=w_0\wedge\ldots\wedge w_k$, the weight of a component $C$ of the cycle $Z(w)$ equals the signed Euler obstruction of the variety $U$ at a generic point of $C$: in particular, if $C\not\subset\sing U$ (which is always the case under $k<\dim\sing U$), then the weight is 1, however otherwise it may be arbitrary. At the same time, this condition of general position on $w$ is far less constructive than the ones elaborated in this text, so it does not look practical.

\vspace{1ex}

2. How to extend our results to the CSM classes $c^{SM}_k(U)$ for $k>\dim\sing U$?

\vspace{1ex}

We do not touch it in this text, though  our affine characteristic classes are designed to study arbitrarily singular varieties (as well as the conic Lagrangian cycles in \cite{maxim}). One may speculate that such an extension would involve tropical characteristic classes $\Trop_k U_\alpha$ of a smooth stratification $U=\bigcup_\alpha U_\alpha$, weighted by the Euler obstructions of $U$ along these strata. 
According to \cite{maxim}, the right object to stduy instead of the CSM classes $c^{SM}_k(U)$ in this generality might be Chern--Mather classes (which equal CSM for $\dim\sing U<k$).

\vspace{1ex}

3. Equally important from our perspective is to constructively understand the multidegrees $v_i$ of a variety $U\subset\CC^n$ (or rather $Z\subset\CC^n$, in the notation and generality of \cite{maxim}). More specifically, $v_i$ can be seen as the degree of the intersection of the likelihood correspondence variety $X\subset\CC^n\times\C^n$ with a generic $i$-plane $L\subset\C^n$. What could be an explicit combinatorial genericity condition for such a rational plane $L$? 

\vspace{1ex}

One could think of some of the fans $\Trop_k U_\alpha$ being immersed or submersed by $L$.

\vspace{1ex}

{\bf Generic torus projections.} Our results apply to a rational form $w:=\frac{dm_0}{m_0}\wedge\cdots\wedge\frac{dm_k}{m_k}$ (and to the corresponding torus projection $p:=(m_0,\ldots,m_k):\CC^n\to\CC^{k+1}$), if the linear projection $dp$ immerses and submerses respective tropical characteristic classes $\Trop_j U$ (or, more generally, $\langle U\rangle_j$ as defined at the beginning of Section \ref{saff}).

\vspace{1ex}

1. Geometry of fibers of generic projections $U|_p$ is addressed in many other results, notably toric and tropical Bertini theorems on irreducibility of such fibers (\cite{fmz}, \cite{tropbert}). To what extent genericity assumptions on the projection $p$ can be explicitly given in a similar way in terms of the fans $\langle U\rangle_j$ for such results?

\vspace{1ex}

2. Funny enough, if $f$ is a generic polynomial with a given Newton polygon, then the planar curve $U:=\{f=0\}\subset\CC^2$ satisfies the results of this paper for any 1-form $w$ (i.e. our genericity condition that the line $\ker w$ intersects the fan $\Trop U=\Trop P$ by one point 0 is redundant). Specifically, for any 1-form $w$, the weighted number of its zeroes on $U$ equals the signed Euler characteristics of $U$. To what extent this observation survives for generic complete intersections or even sch\"on varieties $U\subset\CC^n$ in higher dimensions?

\vspace{1ex}

Of course, the genericity assumptions on the projection are not completely redundant anymore, but one can at least expect to compute how the fan $\Trop Z(w)$ changes as $w$ approaches a forbidden direction.

\vspace{1ex}

{\bf Gauss--Bonnet for non-torus invariant forms.} 

1. Does Theorem \ref{thnoninv} extend to non-exact 1-forms in the following sense? 

\vspace{1ex}

If $A$ satisfies the assumption of the theorem, and $g_1,\ldots,g_n$ are generic linear combinations of its monomials, then the 1-form $w(x):=\sum_i g_i(x)dx_i/x_i$ has $|w|_U=0|=\sum_{k\geq 0}(\Trop A)^k\Trop_k U$ zeroes on $U$. Extending the reasoning after Theorem \ref{thfull} to this setting is non-trivial.

\vspace{1ex}

2. It would be interesting to know the number of zeroes, if $m\in M$ or $A\subset M$ is arbitrary. In case $U$ is a nondegenerate complete intersection, results of \cite{sci} may help here.

\vspace{1ex}

3. In the setting of question 1, it would be interesting to know the number of zeroes, if $g_i$'s are linear combinations of monomials from different generic sets $A_i\subset M$. In case $U$ is a nondegenerate complete intersection, results of \cite{kfmn} may help here.

\vspace{2ex}

{\sc London Institute for Mathematical Sciences, UK}

{\sc Email}: aes@lims.ac.uk

\begin{thebibliography}{99}

\bibitem[AK97]{abrkaru00}
D. Abramovich, K. Karu, Weak semistable reduction in characteristic 0, Invent. Math. 139 (2000) 241--273,  alg-geom/9707012

\bibitem[EGS03]{egs}
W. Ebeling, S. M. Gusein-Zade, J. Seade, Homological index for 1-forms and a Milnor number for isolated singularities, International Journal of Mathematics 15 (2004) 895--905, 	arXiv:math/0307239

\bibitem[ET24]{bd24}
B. El Hilany, E. Tsigaridas, Bounds on the infimum of polynomials over a generic semi-algebraic set using asymptotic critical values (2024) arXiv:2407.17093

\bibitem[E11]{adv}
A. Esterov, The discriminant of a system of equations, Adv. in Math. 245 (2013) 534--572, arXiv:1110.4060

\bibitem[E13]{jems}
A. Esterov, Characteristic classes of affine varieties and Plucker formulas for affine morphisms, J. of the EMS 20 (2018) 15-59, arXiv:1305.3234

\bibitem[E24]{sci}
A. Esterov, Sch\"on complete intersections, arXiv:2401.1209

\bibitem[E06]{kfmn}
A. Esterov, Determinantal singularities and Newton polyhedra, 	arXiv:0906.5097

\bibitem[FMZ14]{fmz}
C. Fuchs, V. Mantova, U. Zannier. On fewnomials, integral points, and a toric version of Bertini’s theorem, J. of the AMS 31 (2018) 107--134, arXiv:1412.4548

\bibitem[FS94]{fs}
W. Fulton, B. Sturmfels, Intersection theory on toric varieties, Topology 36 (1997) 335-353, arXiv:alg-geom/9403002

\bibitem[G+21]{tropbert}
F. Gandini, M. Hering, M. Maclagan, F. Mohammadi, J. Rajchgot, A. Wheeler, J.  Yu,  Toric and tropical Bertini theorems in positive characteristic, In 2021 Fall Southeastern Sectional Meeting (AMS),	arXiv:2111.13214 

\bibitem[GK02]{khgelf}
O. Gelfond, A. Khovanskii, Toric geometry and Grothendieck residues, Mosc. Math. J. 2 (2002) 99--112, {\tt www.math.toronto.edu/askold/grot-res.pdf}

\bibitem[G75]{greuel}
G. M. Greuel, Der Gauss-Manin-zusammenhang isolierter singularit\"aten von
vollst\"andigen durchschnitten, Math. Ann. 214 (1975), 235–266.

\bibitem[HKT07]{fromhubsch}
P. Hacking, S. Keel, J. Tevelev, Stable pair, tropical, and log canonical compact moduli of del Pezzo surfaces, Invent. Math. 178 (2009) 173--227, arXiv:math/0702505

\bibitem[H12]{huh}
J. Huh, The maximum likelihood degree of a very affine variety,
Compositio Math. 149 (2013) 1245-1266, arXiv:1207.0553

\bibitem[H13]{huh2}
J. Huh, Varieties with maximum likelihood degree one, J. of Alg. Statistics 5 (2014) 1-17, arXiv:1301.2732

\bibitem[Kh77]{kh77}
A. G. Khovanskii, Newton polyhedra and the genus of complete intersections, Funct. An. and Its Appl. 12 (1978) 38--46

\bibitem[K04]{valya}
V. Kiritchenko, Chern classes for reductive groups and an adjunction formula, Annales de l’Institut Fourier, 56 (2006) 1225--1256, arXiv:math/0411331

\bibitem[L74]{le}
Le Dung Trang, Computation of the Milnor number of an isolated singularity of a complete intersection, Funkci. An. i Pril. 8 (1974) 45–49

\bibitem[LT81]{lt}
Le Dung Trang, B. Teissier, Vari{\'e}t{\'e}s polaires locales et classes de Chern des vari{\'e}t{\'e}s singulieres, Ann. of Math. 114 (1981) 457--491

\bibitem[MR18]{mr}
G. Mikhalkin, J. Rau, Tropical Geometry, 2018, {\tt math.uniandes.edu.co/$\sim$j.rau/downloads/main.pdf}

\bibitem[MS15]{ms}
D. Maclagan, B. Sturmfels, Introduction to Tropical Geometry, Grad. Studies in Math 161, AMS, 2015

\bibitem[MRWW22]{maxim} L. G. Maxim, J. I. Rodriguez, B. Wang, L. Wu, Logarithmic cotangent bundles, Chern-Mather classes, and the Huh--Sturmfels Involution conjecture, Commun. Pure \& Appl. Math. 77 (2024) 1486-1508, arXiv:2202.00554

\bibitem[P78]{piene}
R. Piene, Polar classes of singular varieties, Annales scientifiques de l'\'Ecole Normale Sup\'erieure, Serie 4 11 (1978) 247--276

\bibitem[ST06]{schuerm}
J. Schürmann, M. Tibar, Index formula for MacPherson cycles of affine algebraic varieties, Tohoku Math. J., Second Series 62 (2010) 29--44, 	arXiv:math/0603338

\bibitem[ST07]{st}
B. Sturmfels, J. Tevelev, Elimination theory for tropical varieties, Math. Res. Lett. 15 (2008) 543--562, arXiv:0704.3471

\bibitem[TW24]{ed24}
S. Telen, M. Wiesmann, Euler Stratifications of Hypersurface Families  (2024) arXiv:2407.18176

\bibitem[T04]{tev}
J. Tevelev, Compactifications of subvarieties of tori, Amer. J. of Math. 129 (2007) 1087--1104, arXiv:math/0412329

\end{thebibliography}
\end{document}